\documentclass{amsart}%
\usepackage{amssymb}
\usepackage{amsfonts}%
\usepackage{amsmath}%
\usepackage{graphicx}
\providecommand{\U}[1]{\protect \rule{.1in}{.1in}}
\newtheorem{theorem}{Theorem}
\theoremstyle{plain}

\newtheorem{corollary}{Corollary}

\newtheorem{definition}{Definition}

\newtheorem{lemma}{Lemma}

\newtheorem{remark}{Remark}

\numberwithin{equation}{section}
\begin{document}
\title[Product generalized local Morrey estimates]{Product generalized local Morrey spaces and commutators of multi-sublinear
operators generated by multilinear {Calder\'{o}n-Zygmund operators} and {local
Campanato functions} }
\author{F. GURBUZ}
\address{ANKARA UNIVERSITY, FACULTY OF SCIENCE, DEPARTMENT OF MATHEMATICS, TANDO\u{G}AN
06100, ANKARA, TURKEY }
\curraddr{}
\email{feritgurbuz84@hotmail.com}
\urladdr{}
\thanks{}
\thanks{}
\thanks{}
\date{}
\dedicatory{ }
\begin{abstract}
The aim of this paper is to get the boundedness of the commutators of
multi-sublinear operators generated by local campanato functions and
multilinear {Calder\'{o}n-Zygmund} operators on the product generalized local
Morrey spaces.

\end{abstract}
\maketitle

\section{Introduction and main results}

Because of the need for study of the local behavior of solutions of second
order elliptic partial differential equations (PDEs) and together with the now
well-studied Sobolev Spaces, constitude a formidable three parameter family of
spaces useful for proving regularity results for solutions to various PDEs,
especially for non-linear elliptic systems, in 1938, Morrey \cite{Morrey}
introduced the classical Morrey spaces $L_{p,\lambda}$ which naturally are
generalizations of the classical Lebesgue spaces.

We will say that a function $f\in L_{p,\lambda}=L_{p,\lambda}\left(
{\mathbb{R}^{n}}\right)  $ if%
\begin{equation}
\sup_{x\in{\mathbb{R}^{n},r>0}}\left[  r^{-\lambda}%
{\displaystyle \int \limits_{B(x,r)}}
\left \vert f\left(  y\right)  \right \vert ^{p}dy\right]  ^{1/p}<\infty
.\label{1.1}%
\end{equation}
Here, $1<p<\infty$ and $0<\lambda<n$ and the quantity of (\ref{1.1}) is the
$\left(  p,\lambda \right)  $-Morrey norm, denoted by $\left \Vert f\right \Vert
_{L_{p,\lambda}}$. We also refer to \cite{Adams} for the latest research on
the theory of Morrey spaces associated with Harmonic Analysis. In recent
years, more and more researches focus on function spaces based on Morrey
spaces to fill in some gaps in the theory of Morrey type spaces (see, for
example, \cite{Gurbuz1, Gurbuz2, Gurbuz3, Gurbuz4, Gurbuz5, Pal}). Moreover,
these spaces are useful in harmonic analysis and PDEs. But, this topic exceeds
the scope of this paper. Thus, we omit the details here.

First of all, we recall some explanations and notations used in the paper.

Recall that the concept of the generalized local (central) Morrey space
$LM_{p,\varphi}^{\{x_{0}\}}$ has been introduced in \cite{BGGS} and studied in
\cite{Gurbuz1, Gurbuz2}.

\begin{definition}
\textbf{(Generalized local (central) Morrey space) }Let $\varphi(x,r)$ be a
positive measurable function on ${\mathbb{R}^{n}}\times(0,\infty)$ and $1\leq
p<\infty$. For any fixed $x_{0}\in{\mathbb{R}^{n}}$ we denote by
$LM_{p,\varphi}^{\{x_{0}\}}\equiv LM_{p,\varphi}^{\{x_{0}\}}({\mathbb{R}^{n}%
})$ the generalized local Morrey space, the space of all functions $f\in
L_{p}^{loc}({\mathbb{R}^{n}})$ with finite quasinorm
\[
\Vert f\Vert_{LM_{p,\varphi}^{\{x_{0}\}}}=\sup \limits_{r>0}\varphi
(x_{0},r)^{-1}\,|B(x_{0},r)|^{-\frac{1}{p}}\, \Vert f\Vert_{L_{p}(B(x_{0}%
,r))}<\infty.
\]

\end{definition}

According to this definition, we recover the local Morrey space $LL_{p,\lambda
}^{\{x_{0}\}}$ under the choice $\varphi(x_{0},r)=r^{\frac{\lambda-n}{p}}$:
\[
LL_{p,\lambda}^{\{x_{0}\}}=LM_{p,\varphi}^{\{x_{0}\}}\mid_{\varphi
(x_{0},r)=r^{\frac{\lambda-n}{p}}}.
\]
For the properties and applications of generalized local (central) Morrey
spaces $LM_{p,\varphi}^{\{x_{0}\}}$, see also \cite{BGGS, Gurbuz1, Gurbuz2}.

On the other hand, in 1976, Coifman et al. \cite{CRW} introduced the
commutator $\overline{T}_{b}$ generated by the Calder\'{o}n-Zygmund operator
$\overline{T}$ and a locally integrable function $b$ as follows:
\begin{equation}
\overline{T}_{b}f(x)=[b,\overline{T}]f\left(  x\right)  \equiv b(x)\overline
{T}f(x)-\overline{T}(bf)(x)=\int \limits_{{\mathbb{R}^{n}}}K\left(  x,y\right)
[b(x)-b(y)]f(y)dy,\label{e3}%
\end{equation}
with the kernel $K$ satisfying the following size condition:%

\[
K\left(  x,y\right)  \leq C\left \vert x-y\right \vert ^{-n},~~~~~~~~\text{
}x\neq y,
\]
and some smoothness assumption. A celebrated result is that $\overline{T}$ is
bounded operator on $L_{p}$ space, where $1<p<\infty$. Sometimes, the
commutator defined by (\ref{e3}) is also called the commutator in Coifman et
al.'s sense, which has its root in the complex analysis and harmonic analysis
(see \cite{CRW}). The main result from \cite{CRW} states that, if and only if
$b\in BMO$ (bounded mean oscillation space), $T_{b}$ is a bounded operator on
$L_{p}\left(  {\mathbb{R}^{n}}\right)  $, $1<p<\infty$. It is worth noting
that for a constant $C$, if $\overline{T}$ is linear we have,%
\begin{align*}
\lbrack b+C,\overline{T}]f  & =\left(  b+C\right)  \overline{T}f-\overline
{T}(\left(  b+C\right)  f)\\
& =b\overline{T}f+C\overline{T}f-\overline{T}\left(  bf\right)  -C\overline
{T}f\\
& =[b,\overline{T}]f.
\end{align*}
This leads one to intuitively look to spaces for which we identify functions
which differ by constants, and so it is no surprise that $b\in BMO$
or{\ }$LC_{q,\lambda}^{\left \{  x_{0}\right \}  }\left(  {\mathbb{R}^{n}%
}\right)  $ (local Campanato space) {has had the most historical significance.
}Also, the definition and some proporties of local Campanato space
$LC_{q,\lambda}^{\left \{  x_{0}\right \}  }\left(  {\mathbb{R}^{n}}\right)  $
that we need in the proof of commutators are as follows.

\begin{definition}
\cite{BGGS, Gurbuz1} Let $1\leq q<\infty$ and $0\leq \lambda<\frac{1}{n}$. A
local Campanato function $b\in L_{q}^{loc}\left(  {\mathbb{R}^{n}}\right)  $
is said to belong to the $LC_{q,\lambda}^{\left \{  x_{0}\right \}  }\left(
{\mathbb{R}^{n}}\right)  $, if%
\begin{equation}
\left \Vert b\right \Vert _{LC_{q,\lambda}^{\left \{  x_{0}\right \}  }}%
=\sup_{r>0}\left(  \frac{1}{\left \vert B\left(  x_{0},r\right)  \right \vert
^{1+\lambda q}}\int \limits_{B\left(  x_{0},r\right)  }\left \vert b\left(
y\right)  -b_{B\left(  x_{0},r\right)  }\right \vert ^{q}dy\right)  ^{\frac
{1}{q}}<\infty,\label{e51}%
\end{equation}
where%
\[
b_{B\left(  x_{0},r\right)  }=\frac{1}{\left \vert B\left(  x_{0},r\right)
\right \vert }\int \limits_{B\left(  x_{0},r\right)  }b\left(  y\right)  dy.
\]
Define%
\[
LC_{q,\lambda}^{\left \{  x_{0}\right \}  }\left(  {\mathbb{R}^{n}}\right)
=\left \{  b\in L_{q}^{loc}\left(  {\mathbb{R}^{n}}\right)  :\left \Vert
b\right \Vert _{LC_{q,\lambda}^{\left \{  x_{0}\right \}  }}<\infty \right \}  .
\]

\end{definition}

\begin{remark}
If two functions which differ by a constant are regarded as a function in the
space $LC_{q,\lambda}^{\left \{  x_{0}\right \}  }\left(  {\mathbb{R}^{n}%
}\right)  $, then $LC_{q,\lambda}^{\left \{  x_{0}\right \}  }\left(
{\mathbb{R}^{n}}\right)  $ becomes a Banach space. The space $LC_{q,\lambda
}^{\left \{  x_{0}\right \}  }\left(  {\mathbb{R}^{n}}\right)  $ when
$\lambda=0$ is just the $LC_{q}^{\left \{  x_{0}\right \}  }({\mathbb{R}^{n}})$.
Apparently, (\ref{e51}) is equivalent to the following condition:%
\[
\sup_{r>0}\inf_{c\in%
\mathbb{C}
}\left(  \frac{1}{\left \vert B\left(  x_{0},r\right)  \right \vert ^{1+\lambda
q}}\int \limits_{B\left(  x_{0},r\right)  }\left \vert b\left(  y\right)
-c\right \vert ^{q}dy\right)  ^{\frac{1}{q}}<\infty.
\]
Also, in \cite{LuYang1}, Lu and Yang have introduced the central BMO space
$CBMO_{q}({\mathbb{R}^{n}})=LC_{q,0}^{\{0\}}({\mathbb{R}^{n}})$. Note that
$BMO({\mathbb{R}^{n}})\subset \bigcap \limits_{q>1}LC_{q}^{\left \{
x_{0}\right \}  }({\mathbb{R}^{n}})$, $1\leq q<\infty$. Moreover, one can
imagine that the behavior of $LC_{q}^{\left \{  x_{0}\right \}  }({\mathbb{R}%
^{n}})$ may be quite different from that of $BMO({\mathbb{R}^{n}})$, since
there is no analogy of the famous John-Nirenberg inequality of
$BMO({\mathbb{R}^{n}})$ for the space $LC_{q}^{\left \{  x_{0}\right \}
}({\mathbb{R}^{n}})$.
\end{remark}

\begin{lemma}
\label{Lemma 4}\cite{BGGS, Gurbuz1} Let $b$ be a local Campanato function in
$LC_{q,\lambda}^{\left \{  x_{0}\right \}  }\left(
\mathbb{R}
^{n}\right)  $, $1\leq q<\infty$, $0\leq \lambda<\frac{1}{n}$ and $r_{1}$,
$r_{2}>0$. Then%
\begin{equation}
\left(  \frac{1}{\left \vert B\left(  x_{0},r_{1}\right)  \right \vert
^{1+\lambda q}}%
{\displaystyle \int \limits_{B\left(  x_{0},r_{1}\right)  }}
\left \vert b\left(  y\right)  -b_{B\left(  x_{0},r_{2}\right)  }\right \vert
^{q}dy\right)  ^{\frac{1}{q}}\leq C\left(  1+\left \vert \ln \frac{r_{1}}{r_{2}%
}\right \vert \right)  \left \Vert b\right \Vert _{LC_{q,\lambda}^{\left \{
x_{0}\right \}  }},\label{a*}%
\end{equation}
where $C>0$ is independent of $b$, $r_{1}$ and $r_{2}$.

From this inequality (\ref{a*}), we have%
\begin{equation}
\left \vert b_{B\left(  x_{0},r_{1}\right)  }-b_{B\left(  x_{0},r_{2}\right)
}\right \vert \leq C\left(  1+\ln \frac{r_{1}}{r_{2}}\right)  \left \vert
B\left(  x_{0},r_{1}\right)  \right \vert ^{\lambda}\left \Vert b\right \Vert
_{LC_{q,\lambda}^{\left \{  x_{0}\right \}  }},\label{b*}%
\end{equation}
and it is easy to see that%
\begin{equation}
\left \Vert b-\left(  b\right)  _{B}\right \Vert _{L_{q}\left(  B\right)  }\leq
C\left(  1+\ln \frac{r_{1}}{r_{2}}\right)  r^{\frac{n}{q}+n\lambda}\left \Vert
b\right \Vert _{LC_{q,\lambda}^{\left \{  x_{0}\right \}  }}.\label{c*}%
\end{equation}

\end{lemma}

\begin{remark}
Let $x_{0}\in{\mathbb{R}^{n}}$, $1<p_{i},q_{i}<\infty$, for $i=1,\ldots,m$
such that $\frac{1}{p}=\frac{1}{p_{1}}+\cdots+\frac{1}{p_{m}}+\frac{1}{q_{1}%
}+\frac{1}{q_{2}}+\cdots+\frac{1}{q_{m}}$ and $\overrightarrow{b}\in
LC_{q_{i},\lambda_{i}}^{\left \{  x_{0}\right \}  }({\mathbb{R}^{n}})$ for
$0\leq \lambda_{i}<\frac{1}{n}$, $i=1,\ldots,m$. Then, from Lemma
\ref{Lemma 4}, it is easy to see that%
\[
\left \Vert b_{i}-\left(  b_{i}\right)  _{B}\right \Vert _{L_{q_{i}}\left(
B\right)  }\leq Cr^{\frac{n}{q_{i}}+n\lambda_{i}}\left \Vert b_{i}\right \Vert
_{LC_{q_{i},\lambda_{i}}^{\left \{  x_{0}\right \}  }},
\]
and%
\begin{equation}
\left \Vert b_{i}-\left(  b_{i}\right)  _{B}\right \Vert _{L_{q_{i}}\left(
2B\right)  }\leq \left \Vert b_{i}-\left(  b_{i}\right)  _{2B}\right \Vert
_{L_{q_{i}}\left(  2B\right)  }+\left \Vert \left(  b_{i}\right)  _{B}-\left(
b_{i}\right)  _{2B}\right \Vert _{L_{q_{i}}\left(  2B\right)  }\lesssim
r^{\frac{n}{q_{i}}+n\lambda_{i}}\left \Vert b_{i}\right \Vert _{LC_{q_{i}%
,\lambda_{i}}^{\left \{  x_{0}\right \}  }},\label{*}%
\end{equation}
for $i=1$, $2$.
\end{remark}

On the other hand, multilinear {Calder\'{o}n-Zygmund theory is a natural
generalization of the linear case. The initial work on the class of
}multilinear {Calder\'{o}n-Zygmund operators has been done by Coifman and
Meyer in \cite{CM}. Moreover, }the study of multilinear singular integrals has
motivated not only as the generalization of the theory of linear ones but also
their natural appearance in analysis. It has received increasing attention and
much development in recent years, such as the study of the bilinear Hilbert
transform by Lacey and Thiele \cite{Lacey 1, Lacey 2} and the systematic
treatment of multilinear {Calder\'{o}n-Zygmund operators by Grafakos-Torres
\cite{Grafakos1, Grafakos2, Grafakos3*} and Grafakos-Kalton \cite{Grafakos3}.
Meanwhile, the commutators generated by the }multilinear singular integral and
$BMO$ functions of Lipschitz functions also attract much attention, since the
commutator is more singular than the singular integral operator itself.

Let ${\mathbb{R}^{n}}$ be the $n$-dimensional Euclidean space of points
$x=(x_{1},...,x_{n})$ with norm $|x|=\left(
{\displaystyle \sum \limits_{i=1}^{n}}
x_{i}^{2}\right)  ^{\frac{1}{2}}$ and $\left(  {\mathbb{R}^{n}}\right)
^{m}={\mathbb{R}^{n}\times \ldots \times \mathbb{R}^{n}}$ be the $m$-fold product
spaces $\left(  m\in%
\mathbb{N}
\right)  $. For $x\in{\mathbb{R}^{n}}$ and $r>0$, we denote by $B(x,r)$ the
open ball centered at $x$ of radius $r$, and by$\,B^{C}(x,r)$ denote its
complement and $|B(x,r)|$ is the Lebesgue measure of the ball $B(x,r)$ and
$|B(x,r)|=v_{n}r^{n}$, where $v_{n}=|B(0,1)|$. Throughout this paper, we
denote by $\overrightarrow{y}=\left(  y_{1},\ldots,y_{m}\right)  $,
$d\overrightarrow{y}=dy_{1}\ldots dy_{m}$, and by $\overrightarrow{f}$ the
$m$-tuple $\left(  f_{1},...,f_{m}\right)  $, $m$, $n$ the nonnegative
integers with $n\geq2$, $m\geq1$.

Suppose that $T^{\left(  m\right)  }$ represents a multilinear or a
multi-sublinear operator, which satisfies that for any $m\in%
\mathbb{N}
$ and $\overrightarrow{f}=\left(  f_{1},\ldots,f_{m}\right)  $, suppose each
$f_{i}$ $\left(  i=1,\ldots,m\right)  $ is integrable on ${\mathbb{R}^{n}}$
with compact support and $x\notin%
{\displaystyle \bigcap \limits_{i=1}^{m}}
suppf_{i}$,%
\begin{equation}
\left \vert T^{\left(  m\right)  }\left(  \overrightarrow{f}\right)  \left(
x\right)  \right \vert \leq c_{0}%
{\displaystyle \int \limits_{\left(  {\mathbb{R}^{n}}\right)  ^{m}}}
\frac{1}{\left \vert \left(  x-y_{1},\ldots,x-y_{m}\right)  \right \vert ^{mn}%
}\left \{
{\displaystyle \prod \limits_{i=1}^{m}}
\left \vert f_{i}\left(  y_{i}\right)  \right \vert \right \}  d\overrightarrow
{y},\label{4}%
\end{equation}
where $c_{0}$ is independent of $\overrightarrow{f}$ and $x$.

We point out that the condition (\ref{4}) in the case of $m=1$ was first
introduced by Soria and Weiss in \cite{SW} . The condition (\ref{4}) is
satisfied by many interesting operators in harmonic analysis, such as the
$m$-linear Calder\'{o}n--Zygmund operators, $m$-sublinear Carleson's maximal
operator, $m$-sublinear Hardy--Littlewood maximal operator, C. Fefferman's
singular multipliers, R. Fefferman's $m$-linear singular integrals,
Ricci--Stein's $m$-linear oscillatory singular integrals, the $m$-linear
Bochner--Riesz means and so on (see \cite{BGGS, Gurbuz1, Gurbuz2, SW} for details).

We are going to be working on ${\mathbb{R}^{n}}$. Let's begin with the
recalling of the multilinear {Calder\'{o}n-Zygmund operator }$\overline
{T}^{\left(  m\right)  }$ $\left(  m\in%
\mathbb{N}
\right)  ${. }Let $\overrightarrow{f}\in L_{p_{1}}^{loc}\left(  {\mathbb{R}%
^{n}}\right)  \times \ldots \times L_{p_{m}}^{loc}\left(  {\mathbb{R}^{n}%
}\right)  $. {The }$m$(multi)-linear {Calder\'{o}n-Zygmund }operator
$\overline{T}^{\left(  m\right)  }$ is defined by
\[
\overline{T}^{\left(  m\right)  }\left(  \overrightarrow{f}\right)  \left(
x\right)  =\overline{T}^{\left(  m\right)  }\left(  f_{1},\ldots,f_{m}\right)
\left(  x\right)  =%
{\displaystyle \int \limits_{\left(  {\mathbb{R}^{n}}\right)  ^{m}}}
K\left(  x,y_{1},\ldots,y_{m}\right)  \left \{
{\displaystyle \prod \limits_{i=1}^{m}}
f_{i}\left(  y_{i}\right)  \right \}  dy_{1}\cdots dy_{m},
\]
for test vector $\overrightarrow{f}=\left(  f_{1},\ldots,f_{m}\right)  $ and
$x\notin%
{\displaystyle \bigcap \limits_{i=1}^{m}}
suppf_{i}$, where $K$ is an $m$-{Calder\'{o}n-Zygmund kernel which is a
locally integrable function defined }away from the diagonal $y_{0}%
=y_{1}=\cdots=y_{m}$ on $\left(  {\mathbb{R}^{n}}\right)  ^{m+1}$ satisfying
the {following size estimate:}%
\begin{equation}
\left \vert K\left(  x,y_{1},\ldots,y_{m}\right)  \right \vert \leq \frac
{C}{\left \vert \left(  x-y_{1},\ldots,x-y_{m}\right)  \right \vert ^{mn}%
},\label{10*}%
\end{equation}
for some $C>0$ and some smoothness estimates, see \cite{Grafakos3}%
-\cite{Grafakos1} for details.

At the same time, Grafakos and Torres \cite{Grafakos1, Grafakos3*} have proved
that the multilinear {Calder\'{o}n-Zygmund operator is bounded on the product
of Lebesgue spaces. }

\begin{theorem}
\label{teo0}\cite{Grafakos1, Grafakos3*} Let $\overline{T}^{\left(  m\right)
} $ is a $m$-linear {Calder\'{o}n-Zygmund operator. Then, for any numbers
}$1\leq p_{1},\ldots,p_{m}<\infty$ with $\frac{1}{p}=\frac{1}{p_{1}}%
+\cdots+\frac{1}{p_{m}}$, $\overline{T}^{\left(  m\right)  }$ can be extended
to a bounded operator from $L_{p_{1}}\times \cdots \times L_{p_{m}}$ into
$L_{p}$, and bounded from $L_{1}\times \cdots \times L_{1}$ into $L_{\frac{1}%
{m},\infty}$.
\end{theorem}

Recently, Xu \cite{Xu} has established the boundedness on the{\ product of
Lebesgue space} for the commutators generated by $m$-linear
{Calder\'{o}n-Zygmund }singular integrals and $RBMO$ functions with
nonhomogeneous. Inspired by \cite{Grafakos1}, \cite{Grafakos3*}, \cite{Xu}, we
will introduce the commutators $\overline{T}_{\overrightarrow{b}}^{\left(
m\right)  }$ generated by $m$-linear {Calder\'{o}n-Zygmund operators
}$\overline{T}^{\left(  m\right)  }$ {and local Campanato functions
}$\overrightarrow{b}=\left(  b_{1},\ldots,b_{m}\right)  ${\ }%
\[
\overline{T}_{\overrightarrow{b}}^{\left(  m\right)  }\left(  \overrightarrow
{f}\right)  \left(  x\right)  =%
{\displaystyle \int \limits_{\left(  {\mathbb{R}^{n}}\right)  ^{m}}}
K\left(  x,y_{1},\ldots,y_{m}\right)  \left[
{\displaystyle \prod \limits_{i=1}^{m}}
\left[  b_{i}\left(  x\right)  -b_{i}\left(  y_{i}\right)  \right]
f_{i}\left(  y_{i}\right)  \right]  d\overrightarrow{y},
\]
where $K\left(  x,y_{1},\ldots,y_{m}\right)  $ is a $m$-linear
{Calder\'{o}n-Zygmund kernel, }$b_{i}\in LC_{q_{i},\lambda_{i}}^{\left \{
x_{0}\right \}  }({\mathbb{R}^{n}})$ (local Campanato spaces) for $0\leq
\lambda_{i}<\frac{1}{n}$, $i=1,\ldots,m$. We would like to point out that
$\overline{T}_{b}$ is the special case of $\overline{T}_{\overrightarrow{b}%
}^{\left(  m\right)  }$ with taking $m=1$.

Closely related to the above results, in this paper in the case of $b_{i}\in
LC_{q_{i},\lambda_{i}}^{\left \{  x_{0}\right \}  }({\mathbb{R}^{n}})$ for
$0\leq \lambda_{i}<\frac{1}{n}$, $i=1,\ldots,m$, we find the sufficient
conditions on $(\varphi_{1},\ldots \varphi_{m},\varphi)$ which ensures the
boundedness of the commutator operators $T_{\overrightarrow{b}}^{\left(
m\right)  }$ from $LM_{p_{1},\varphi_{1}}^{\{x_{0}\}}\times \cdots$ $\times$
$LM_{p_{m},\varphi_{m}}^{\{x_{0}\}}$ to $LM_{p,\varphi}^{\{x_{0}\}}$, where
$1<p_{i},q_{i}<\infty,$ for $i=1,\ldots,m$ such that $\frac{1}{p}=\frac
{1}{p_{1}}+\cdots+\frac{1}{p_{m}}+\frac{1}{q_{1}}+\frac{1}{q_{2}}+\cdots
+\frac{1}{q_{m}}$. In fact, in this paper the results of \cite{BGGS} and
\cite{Gurbuz1} (by taking $\Omega \equiv1$ there) will be generalized to the
multilinear case; we omit the details here.

\begin{remark}
Our results in this paper remain true for the inhomogeneous versions of local
Campanato spaces $LC_{q,\lambda}^{\left \{  x_{0}\right \}  }\left(
{\mathbb{R}^{n}}\right)  $ for $0\leq \lambda<\frac{1}{n}$ and generalized
local Morrey spaces $LM_{p,\varphi}^{\{x_{0}\}}$.
\end{remark}

We now make some conventions. Throughout this paper, we use the symbol
$A\lesssim B$ to denote that there exists a positive consant $C$ such that
$A\leq CB$. If $A\lesssim B$ and $B\lesssim A$, we then write $A\approx B$ and
say that $A$ and $B$ are equivalent. For a fixed $p\in \left[  1,\infty \right)
$, $p^{\prime}$ denotes the dual or conjugate exponent of $p$, namely,
$p^{\prime}=\frac{p}{p-1}$ and we use the convention $1^{\prime}=\infty$ and
$\infty^{\prime}=1$.

Our main results can be formulated as follows.

\begin{theorem}
\label{Teo 5}Let $x_{0}\in{\mathbb{R}^{n}}$, $1<p_{i},q_{i}<\infty$, for
$i=1,\ldots,m$ such that $\frac{1}{p}=\frac{1}{p_{1}}+\cdots+\frac{1}{p_{m}%
}+\frac{1}{q_{1}}+\frac{1}{q_{2}}+\cdots+\frac{1}{q_{m}}$ and $\overrightarrow
{b}\in LC_{q_{i},\lambda_{i}}^{\left \{  x_{0}\right \}  }({\mathbb{R}^{n}})$
for $0\leq \lambda_{i}<\frac{1}{n}$, $i=1,\ldots,m$. Let also, $T^{\left(
m\right)  }$ $\left(  m\in%
\mathbb{N}
\right)  $ be a multilinear operator satisfying condition (\ref{4}), bounded
from $L_{p_{1}}\times \cdots \times L_{p_{m}}$ into $L_{p}$ for $p_{i}>1\left(
i=1,\ldots,m\right)  $. Then the inequality
\begin{equation}
\Vert T_{\overrightarrow{b}}^{\left(  m\right)  }\left(  \overrightarrow
{f}\right)  \Vert_{L_{p}(B(x_{0},r))}\lesssim%
{\displaystyle \prod \limits_{i=1}^{m}}
\Vert \overrightarrow{b}\Vert_{LC_{q_{i},\lambda_{i}}^{\left \{  x_{0}\right \}
}}\,r^{\frac{n}{p}}\int \limits_{2r}^{\infty}\left(  1+\ln \frac{t}{r}\right)
^{m}%
{\displaystyle \prod \limits_{i=1}^{m}}
\Vert f_{i}\Vert_{L_{p_{i}}(B(x_{0},t))}\frac{dt}{t^{^{n\left(
{\displaystyle \sum \limits_{i=1}^{n}}
\frac{1}{p_{i}}-%
{\displaystyle \sum \limits_{i=1}^{n}}
\lambda_{i}\right)  +1}}}\label{200}%
\end{equation}
holds for any ball $B(x_{0},r)$ and for all $\overrightarrow{f}\in L_{p_{1}%
}^{loc}\left(  {\mathbb{R}^{n}}\right)  \times \cdots L_{p_{m}}^{loc}\left(
{\mathbb{R}^{n}}\right)  $.
\end{theorem}

\begin{theorem}
\label{teo15}Let $x_{0}\in{\mathbb{R}^{n}}$, $1<p_{i},q_{i}<\infty$, for
$i=1,\ldots,m$ such that $\frac{1}{p}=\frac{1}{p_{1}}+\cdots+\frac{1}{p_{m}%
}+\frac{1}{q_{1}}+\frac{1}{q_{2}}+\cdots+\frac{1}{q_{m}}$ and $\overrightarrow
{b}\in LC_{q_{i},\lambda_{i}}^{\left \{  x_{0}\right \}  }({\mathbb{R}^{n}})$
for $0\leq \lambda_{i}<\frac{1}{n}$, $i=1,\ldots,m$. Let also, $T^{\left(
m\right)  }$ $\left(  m\in%
\mathbb{N}
\right)  $ be a multilinear operator satisfying condition (\ref{4}), bounded
from $L_{p_{1}}\times \cdots \times L_{p_{m}}$ into $L_{p}$ for $p_{i}>1\left(
i=1,\ldots,m\right)  $. If functions $\varphi,\varphi_{i}:{\mathbb{R}%
^{n}\times}\left(  0,\infty \right)  \rightarrow \left(  0,\infty \right)  ,$
$\left(  i=1,\ldots,m\right)  $ and $(\varphi_{1},\ldots \varphi_{m},\varphi)$
satisfies the condition%
\begin{equation}%
{\displaystyle \int \limits_{r}^{\infty}}
\left(  1+\ln \frac{t}{r}\right)  ^{m}\frac{\operatorname*{essinf}%
\limits_{t<\tau<\infty}%
{\displaystyle \prod \limits_{i=1}^{m}}
\varphi_{i}(x_{0},\tau)\tau^{\frac{n}{p_{i}}}}{t^{n\left(
{\displaystyle \sum \limits_{i=1}^{n}}
\frac{1}{p_{i}}-%
{\displaystyle \sum \limits_{i=1}^{n}}
\lambda_{i}\right)  +1}}dt\leq C\varphi \left(  x_{0},r\right)  ,\label{50}%
\end{equation}
where $C$ does not depend on $r$.

Then the operator $T_{\overrightarrow{b}}^{\left(  m\right)  }$ $\left(  m\in%
\mathbb{N}
\right)  $ is bounded from product space $LM_{p_{1},\varphi_{1}}^{\{x_{0}%
\}}\times \cdots$ $\times$ $LM_{p_{m},\varphi_{m}}^{\{x_{0}\}}$ to
$LM_{p,\varphi}^{\{x_{0}\}}$ for $p_{i}>1\left(  i=1,\ldots,m\right)  $.
Moreover, we have for $p_{i}>1\left(  i=1,\ldots,m\right)  $
\begin{equation}
\left \Vert T_{\overrightarrow{b}}^{\left(  m\right)  }\left(  \overrightarrow
{f}\right)  \right \Vert _{LM_{p,\varphi}^{\{x_{0}\}}}\lesssim%
{\displaystyle \prod \limits_{i=1}^{m}}
\left \Vert \overrightarrow{b}\right \Vert _{LC_{q_{i},\lambda_{i}}^{\left \{
x_{0}\right \}  }}%
{\displaystyle \prod \limits_{i=1}^{m}}
\left \Vert f_{i}\right \Vert _{LM_{p_{i},\varphi_{i}}^{\{x_{0}\}}}.\label{51}%
\end{equation}

\end{theorem}

For the $m$-sublinear commutator of the $m$-sublinear maximal operator%

\[
M_{\overrightarrow{b}}^{\left(  m\right)  }\left(  \overrightarrow{f}\right)
\left(  x\right)  =\sup_{t>0}\frac{1}{\left \vert B\left(  x,t\right)
\right \vert }\int \limits_{B\left(  x,t\right)  }%
{\displaystyle \prod \limits_{i=1}^{m}}
\left[  b_{i}\left(  x\right)  -b_{i}\left(  y_{i}\right)  \right]  \left \vert
f_{i}\left(  y_{i}\right)  \right \vert d\overrightarrow{y}%
\]

from Theorem \ref{teo15} we get the following new results.

\begin{corollary}
\label{corollary 3*}Let $x_{0}\in{\mathbb{R}^{n}}$, $1<p_{i},q_{i}<\infty$,
for $i=1,\ldots,m$ such that $\frac{1}{p}=\frac{1}{p_{1}}+\cdots+\frac
{1}{p_{m}}+\frac{1}{q_{1}}+\frac{1}{q_{2}}+\cdots+\frac{1}{q_{m}}$ and
$\overrightarrow{b}\in LC_{q_{i},\lambda_{i}}^{\left \{  x_{0}\right \}
}({\mathbb{R}^{n}})$ for $0\leq \lambda_{i}<\frac{1}{n}$, $i=1,\ldots,m$ and
$(\varphi_{1},\ldots \varphi_{m},\varphi)$ satisfies condition (\ref{50}).
Then, the operators $M_{\overrightarrow{b}}^{\left(  m\right)  }$ and
$\overline{T}_{\overrightarrow{b}}^{\left(  m\right)  }$ $\left(  m\in%
\mathbb{N}
\right)  $ are bounded from product space $LM_{p_{1},\varphi_{1}}^{\{x_{0}%
\}}\times \cdots$ $\times$ $LM_{p_{m},\varphi_{m}}^{\{x_{0}\}}$ to
$LM_{p,\varphi}^{\{x_{0}\}}$ for $p_{i}>1\left(  i=1,\ldots,m\right)  $.
\end{corollary}

\begin{remark}
Note that, in the case of $m=1$ Theorem \ref{teo15} and Corollary
\ref{corollary 3*} have been proved in \cite{BGGS, Gurbuz1}.
\end{remark}

\begin{corollary}
\label{corollary 10}Let $1<p_{i},q_{i}<\infty$, for $i=1,\ldots,m$ such that
$\frac{1}{p}=\frac{1}{p_{1}}+\cdots+\frac{1}{p_{m}}$ and $\overrightarrow
{b}\in BMO^{m}({\mathbb{R}^{n}})$ for $i=1,\ldots,m$. Let also, $T^{\left(
m\right)  }$ $\left(  m\in%
\mathbb{N}
\right)  $ be a multilinear operator satisfying condition (\ref{4}), bounded
from $L_{p_{1}}\times \cdots \times L_{p_{m}}$ into $L_{p}$ for $p_{i}>1\left(
i=1,\ldots,m\right)  $. If functions $\varphi,\varphi_{i}:{\mathbb{R}%
^{n}\times}\left(  0,\infty \right)  \rightarrow \left(  0,\infty \right)  ,$
$\left(  i=1,\ldots,m\right)  $ and $(\varphi_{1},\ldots \varphi_{m},\varphi)$
satisfies the condition%
\begin{equation}%
{\displaystyle \int \limits_{r}^{\infty}}
\left(  1+\ln \frac{t}{r}\right)  ^{m}\frac{\operatorname*{essinf}%
\limits_{t<\tau<\infty}%
{\displaystyle \prod \limits_{i=1}^{m}}
\varphi_{i}(x,\tau)\tau^{\frac{n}{p_{i}}}}{t^{n%
{\displaystyle \sum \limits_{i=1}^{n}}
\frac{1}{p_{i}}+1}}dt\leq C\varphi \left(  x,r\right)  ,\label{e57}%
\end{equation}
where $C$ does not depend on $r$.

Then the operator $T_{\overrightarrow{b}}^{\left(  m\right)  }$ $\left(  m\in%
\mathbb{N}
\right)  $ is bounded from product space $M_{p_{1},\varphi_{1}}\times \cdots$
$\times$ $M_{p_{m},\varphi_{m}}$ to $M_{p,\varphi}$ for $p_{i}>1\left(
i=1,\ldots,m\right)  $. Moreover, we have for $p_{i}>1\left(  i=1,\ldots
,m\right)  $
\[
\left \Vert T_{\overrightarrow{b}}^{\left(  m\right)  }\left(  \overrightarrow
{f}\right)  \right \Vert _{M_{p,\varphi}}\lesssim%
{\displaystyle \prod \limits_{i=1}^{m}}
\left \Vert \overrightarrow{b}\right \Vert _{BMO^{m}}%
{\displaystyle \prod \limits_{i=1}^{m}}
\left \Vert f_{i}\right \Vert _{M_{p_{i},\varphi_{i}}}.
\]

\end{corollary}

\begin{corollary}
\label{corollary 13} Let $1<p_{i},q_{i}<\infty$, for $i=1,\ldots,m$ such that
$\frac{1}{p}=\frac{1}{p_{1}}+\cdots+\frac{1}{p_{m}}$ and $\overrightarrow
{b}\in BMO^{m}({\mathbb{R}^{n}})$ for $i=1,\ldots,m$ and also $(\varphi
_{1},\ldots \varphi_{m},\varphi)$ satisfies condition (\ref{e57}). Then, the
operators $M_{\overrightarrow{b}}^{\left(  m\right)  }$ and $\overline
{T}_{\overrightarrow{b}}^{\left(  m\right)  }$ $\left(  m\in%
\mathbb{N}
\right)  $ are bounded from product space $M_{p_{1},\varphi_{1}}\times \cdots$
$\times$ $M_{p_{m},\varphi_{m}}$ to $M_{p,\varphi}$ for $p_{i}>1\left(
i=1,\ldots,m\right)  $.
\end{corollary}

\section{Proofs of the main results}

\subsection{\textbf{Proof of Theorem \ref{Teo 5}.}}

\begin{proof}
In order to simplify the proof, we consider only the situation when $m=2$.
Actually, a similar procedure works for all $m\in%
\mathbb{N}
$. Thus, without loss of generality, it is suffice to show that the conclusion
holds for $T_{\overrightarrow{b}}^{\left(  2\right)  }\left(  \overrightarrow
{f}\right)  =T_{\left(  b_{1},b_{2}\right)  }^{\left(  2\right)  }\left(
f_{1},f_{2}\right)  $.

We just consider the case $p_{i},q_{i}>1$ for $i=1,2$. For any $x_{0}%
\in{\mathbb{R}^{n}}$, set $B=B\left(  x_{0},r\right)  $ for the ball centered
at $x_{0}$ and of radius $r$ and $2B=B\left(  x_{0},2r\right)  $. Thus, we
have the following decomposition,%
\begin{align*}
T_{\left(  b_{1},b_{2}\right)  }^{\left(  2\right)  }\left(  f_{1}%
,f_{2}\right)  \left(  x\right)   & =\left[  \left(  b_{1}\left(  x\right)
-\left \{  b_{1}\right \}  _{B}\right)  \right]  \left[  \left(  b_{2}\left(
x\right)  -\left \{  b_{2}\right \}  _{B}\right)  \right]  T^{\left(  2\right)
}\left(  f_{1},f_{2}\right)  \left(  x\right) \\
& -\left[  \left(  b_{1}\left(  x\right)  -\left \{  b_{1}\right \}
_{B}\right)  \right]  T^{\left(  2\right)  }\left[  f_{1},\left(  b_{2}\left(
\cdot \right)  -\left \{  b_{2}\right \}  _{B}\right)  f_{2}\right]  \left(
x\right) \\
& -\left[  \left(  b_{2}\left(  x\right)  -\left \{  b_{2}\right \}
_{B}\right)  \right]  T^{\left(  2\right)  }\left[  \left(  b_{1}\left(
\cdot \right)  -\left \{  b_{1}\right \}  _{B}\right)  f_{1},f_{2}\right]
\left(  x\right) \\
& +T^{\left(  2\right)  }\left[  \left(  b_{1}\left(  \cdot \right)  -\left \{
b_{1}\right \}  _{B}\right)  f_{1},\left(  b_{2}\left(  \cdot \right)  -\left \{
b_{2}\right \}  _{B}\right)  f_{2}\right]  \left(  x\right) \\
& \equiv H_{1}\left(  x\right)  +H_{2}\left(  x\right)  +H_{3}\left(
x\right)  +H_{4}\left(  x\right)  .
\end{align*}

Thus,%
\begin{equation}
\left \Vert T_{\left(  b_{1},b_{2}\right)  }^{\left(  2\right)  }\left(
f_{1},f_{2}\right)  \right \Vert _{L_{p}\left(  B\left(  x_{0},r\right)
\right)  }=\left(
{\displaystyle \int \limits_{B}}
\left \vert T_{\left(  b_{1},b_{2}\right)  }^{\left(  2\right)  }\left(
f_{1},f_{2}\right)  \left(  x\right)  \right \vert ^{p}dx\right)  ^{\frac{1}%
{p}}\leq%
{\displaystyle \sum \limits_{i=1}^{4}}
\left(
{\displaystyle \int \limits_{B}}
\left \vert H_{i}\left(  x\right)  \right \vert ^{p}dx\right)  ^{\frac{1}{p}}=%
{\displaystyle \sum \limits_{i=1}^{4}}
G_{i}.\label{100}%
\end{equation}

One observes that the estimate of $G_{2}$ is analogous to that of $G_{3}$.
Thus, we will only estimate $G_{1}$, $G_{2}$ and $G_{4}$.

Indeed, we also decompose $f_{i}$ as $f_{i}\left(  y_{i}\right)  =$
$f_{i}\left(  y_{i}\right)  \chi_{2B}+f_{i}\left(  y_{i}\right)  \chi_{\left(
2B\right)  ^{C}}$ for $i=1,2$. And, we write $f_{1}=f_{1}^{0}+f_{1}^{\infty}$
and $f_{2}=f_{2}^{0}+f_{2}^{\infty}$, where \ $f_{i}^{0}=f_{i}\chi_{2B}$,
\ $f_{i}^{\infty}=f_{i}\chi_{\left(  2B\right)  ^{C}}$, for $i=1,2$.

$\left(  i\right)  $ For $G_{1}=\left \Vert T_{\left(  b_{1},b_{2}\right)
}^{\left(  2\right)  }\left(  f_{1}^{0},f_{2}^{0}\right)  \right \Vert
_{L_{p}\left(  B\left(  x_{0},r\right)  \right)  }$, we decompose it into four
parts as follows:%
\begin{align*}
G_{1}  & \lesssim \left \Vert \left[  \left(  b_{1}-\left \{  b_{1}\right \}
_{B}\right)  \right]  \left[  \left(  b_{2}-\left \{  b_{2}\right \}
_{B}\right)  \right]  T^{\left(  2\right)  }\left(  f_{1}^{0},f_{2}%
^{0}\right)  \right \Vert _{L_{p}\left(  B\left(  x_{0},r\right)  \right)  }\\
& +\left \Vert \left[  \left(  b_{1}-\left \{  b_{1}\right \}  _{B}\right)
\right]  T^{\left(  2\right)  }\left[  f_{1}^{0},\left(  b_{2}-\left \{
b_{2}\right \}  _{B}\right)  f_{2}^{0}\right]  \right \Vert _{L_{p}\left(
B\left(  x_{0},r\right)  \right)  }\\
& +\left \Vert \left[  \left(  b_{2}-\left \{  b_{2}\right \}  _{B}\right)
\right]  T^{\left(  2\right)  }\left[  \left(  b_{1}-\left \{  b_{1}\right \}
_{B}\right)  f_{1}^{0},f_{2}^{0}\right]  \right \Vert _{L_{p}\left(  B\left(
x_{0},r\right)  \right)  }\\
& +\left \Vert T^{\left(  2\right)  }\left[  \left(  b_{1}-\left \{
b_{1}\right \}  _{B}\right)  f_{1}^{0},\left(  b_{2}-\left \{  b_{2}\right \}
_{B}\right)  f_{2}^{0}\right]  \right \Vert _{L_{p}\left(  B\left(
x_{0},r\right)  \right)  }\\
& \equiv G_{11}+G_{12}+G_{13}+G_{14}.
\end{align*}

Firstly, $1<\overline{p},\overline{q}<\infty$, such that $\frac{1}%
{\overline{p}}=\frac{1}{p_{1}}+\frac{1}{p_{2}}$ and $\frac{1}{\overline{q}%
}=\frac{1}{q_{1}}+\frac{1}{q_{2}}$. Then, using H\"{o}lder's inequality and
from the boundedness of $T^{\left(  2\right)  }$ from $L_{p_{1}}\times
L_{p_{2}}$ into $L_{\overline{p}}$ (see Theorem \ref{teo0}) it follows that:%
\begin{align*}
G_{11}  & \lesssim \left \Vert \left[  \left(  b_{1}-\left \{  b_{1}\right \}
_{B}\right)  \right]  \left[  \left(  b_{2}-\left \{  b_{2}\right \}
_{B}\right)  \right]  \right \Vert _{L_{\overline{q}}\left(  B\right)
}\left \Vert T^{\left(  2\right)  }\left(  f_{1}^{0},f_{2}^{0}\right)
\right \Vert _{L_{\overline{p}}\left(  B\right)  }\\
& \lesssim \left \Vert \left(  b_{1}-\left \{  b_{1}\right \}  _{B}\right)
\right \Vert _{L_{q_{1}}\left(  B\right)  }\left \Vert \left(  b_{2}-\left \{
b_{2}\right \}  _{B}\right)  \right \Vert _{L_{q_{2}}\left(  B\right)
}\left \Vert f_{1}\right \Vert _{L_{p_{1}}\left(  2B\right)  }\left \Vert
f_{2}\right \Vert _{L_{p_{2}}\left(  2B\right)  }\\
& \lesssim \left \Vert \left(  b_{1}-\left \{  b_{1}\right \}  _{B}\right)
\right \Vert _{L_{q_{1}}\left(  B\right)  }\left \Vert \left(  b_{2}-\left \{
b_{2}\right \}  _{B}\right)  \right \Vert _{L_{q_{2}}\left(  B\right)
}r^{n\left(  \frac{1}{p_{1}}+\frac{1}{p_{2}}\right)  }\int \limits_{2r}%
^{\infty}%
{\displaystyle \prod \limits_{i=1}^{2}}
\Vert f_{i}\Vert_{L_{p_{i}}(B(x_{0},t))}\frac{dt}{t^{n\left(  \frac{1}{p_{1}%
}+\frac{1}{p_{2}}\right)  +1}}\\
& \lesssim%
{\displaystyle \prod \limits_{i=1}^{2}}
\Vert \overrightarrow{b}\Vert_{LC_{q_{i},\lambda_{i}}^{\left \{  x_{0}\right \}
}}r^{\frac{n}{p}}\int \limits_{2r}^{\infty}\left(  1+\ln \frac{t}{r}\right)
^{2}%
{\displaystyle \prod \limits_{i=1}^{2}}
\Vert f_{i}\Vert_{L_{p_{i}}(B(x_{0},t))}\frac{dt}{t^{n\left(  \left(  \frac
{1}{p_{1}}+\frac{1}{p_{2}}\right)  -\left(  \lambda_{1}+\lambda_{2}\right)
\right)  +1}}.
\end{align*}
Secondly, for $G_{12}$, let $1<\tau<\infty$, such that $\frac{1}{p}=\frac
{1}{q_{1}}+\frac{1}{\tau}$. Then similar to the estimates for $G_{11}$, we
have%
\begin{align*}
G_{12}  & \lesssim \left \Vert \left(  b_{1}-\left \{  b_{1}\right \}
_{B}\right)  \right \Vert _{L_{q_{1}}\left(  B\right)  }\left \Vert T^{\left(
2\right)  }\left[  f_{1}^{0},\left(  b_{2}-\left \{  b_{2}\right \}
_{B}\right)  f_{2}^{0}\right]  \right \Vert _{L_{\tau}\left(  B\right)  }\\
& \lesssim \left \Vert \left(  b_{1}-\left \{  b_{1}\right \}  _{B}\right)
\right \Vert _{L_{q_{1}}\left(  B\right)  }\left \Vert f_{1}^{0}\right \Vert
_{_{L_{p_{1}}\left(
\mathbb{R}
^{n}\right)  }}\left \Vert \left(  b_{2}-\left \{  b_{2}\right \}  _{2B}\right)
f_{2}^{0}\right \Vert _{L_{k}\left(
\mathbb{R}
^{n}\right)  }\\
& \lesssim \left \Vert \left(  b_{1}-\left \{  b_{1}\right \}  _{B}\right)
\right \Vert _{L_{q_{1}}\left(  B\right)  }\left \Vert \left(  b_{2}-\left \{
b_{2}\right \}  _{B}\right)  \right \Vert _{L_{q_{2}}\left(  2B\right)
}\left \Vert f_{1}\right \Vert _{L_{p_{1}}\left(  2B\right)  }\left \Vert
f_{2}\right \Vert _{L_{p_{2}}\left(  2B\right)  },
\end{align*}
where $1<k<\infty$, such that $\frac{1}{k}=\frac{1}{p_{2}}+\frac{1}{q_{2}%
}=\frac{1}{\tau}-\frac{1}{p_{1}}$.

Hence, we get%
\[
G_{12}\lesssim%
{\displaystyle \prod \limits_{i=1}^{2}}
\Vert \overrightarrow{b}\Vert_{LC_{q_{i},\lambda_{i}}^{\left \{  x_{0}\right \}
}}r^{\frac{n}{p}}\int \limits_{2r}^{\infty}\left(  1+\ln \frac{t}{r}\right)
^{2}%
{\displaystyle \prod \limits_{i=1}^{2}}
\Vert f_{i}\Vert_{L_{p_{i}}(B(x_{0},t))}\frac{dt}{t^{n\left(  \left(  \frac
{1}{p_{1}}+\frac{1}{p_{2}}\right)  -\left(  \lambda_{1}+\lambda_{2}\right)
\right)  +1}}.
\]

Similarly, $G_{13}$ has the same estimate above, here we omit the details,
thus following inequality%
\[
G_{13}\lesssim%
{\displaystyle \prod \limits_{i=1}^{2}}
\Vert \overrightarrow{b}\Vert_{LC_{q_{i},\lambda_{i}}^{\left \{  x_{0}\right \}
}}r^{\frac{n}{p}}\int \limits_{2r}^{\infty}\left(  1+\ln \frac{t}{r}\right)
^{2}%
{\displaystyle \prod \limits_{i=1}^{2}}
\Vert f_{i}\Vert_{L_{p_{i}}(B(x_{0},t))}\frac{dt}{t^{n\left(  \left(  \frac
{1}{p_{1}}+\frac{1}{p_{2}}\right)  -\left(  \lambda_{1}+\lambda_{2}\right)
\right)  +1}}%
\]
is valid.

At last, we consider the term $G_{14}$. Let $1<\tau_{1},\tau_{2}<\infty$, such
that $\frac{1}{\tau_{1}}=\frac{1}{p_{1}}+\frac{1}{q_{1}}$ and $\frac{1}%
{\tau_{2}}=\frac{1}{p_{2}}+\frac{1}{q_{2}}$. It is easy to see that $\frac
{1}{p}=\frac{1}{\tau_{1}}+\frac{1}{\tau_{2}}$. Then by the boundedness of
$T^{\left(  2\right)  }$ from $L_{\tau_{1}}\times L_{\tau_{2}}$ into $L_{p}$
(see Theorem \ref{teo0}), H\"{o}lder's inequality and (\ref{*}), we obtain%
\begin{align*}
G_{14}  & \lesssim \left \Vert \left(  b_{1}-\left \{  b_{1}\right \}
_{B}\right)  f_{1}^{0}\right \Vert _{L_{\tau_{1}}\left(  {\mathbb{R}^{n}%
}\right)  }\left \Vert \left(  b_{2}-\left \{  b_{2}\right \}  _{B}\right)
f_{2}^{0}\right \Vert _{L_{\tau_{2}}\left(  {\mathbb{R}^{n}}\right)  }\\
& \lesssim \left \Vert \left(  b_{1}-\left \{  b_{1}\right \}  _{B}\right)
\right \Vert _{L_{q_{1}}\left(  2B\right)  }\left \Vert \left(  b_{2}-\left \{
b_{2}\right \}  _{B}\right)  \right \Vert _{L_{q_{2}}\left(  2B\right)
}\left \Vert f_{1}\right \Vert _{L_{p_{1}}\left(  2B\right)  }\left \Vert
f_{2}\right \Vert _{L_{p_{2}}\left(  2B\right)  }\\
& \lesssim%
{\displaystyle \prod \limits_{i=1}^{2}}
\Vert \overrightarrow{b}\Vert_{LC_{q_{i},\lambda_{i}}^{\left \{  x_{0}\right \}
}}r^{\frac{n}{p}}\int \limits_{2r}^{\infty}\left(  1+\ln \frac{t}{r}\right)
^{2}%
{\displaystyle \prod \limits_{i=1}^{2}}
\Vert f_{i}\Vert_{L_{p_{i}}(B(x_{0},t))}\frac{dt}{t^{n\left(  \left(  \frac
{1}{p_{1}}+\frac{1}{p_{2}}\right)  -\left(  \lambda_{1}+\lambda_{2}\right)
\right)  +1}}.
\end{align*}
Combining all the estimates of $G_{11}$, $G_{12}$, $G_{13}$, $G_{14}$; there
is%
\begin{align*}
G_{1}  & =\left \Vert T_{\left(  b_{1},b_{2}\right)  }^{\left(  2\right)
}\left(  f_{1}^{0},f_{2}^{0}\right)  \right \Vert _{L_{p}\left(  B\left(
x_{0},r\right)  \right)  }\lesssim%
{\displaystyle \prod \limits_{i=1}^{2}}
\Vert \overrightarrow{b}\Vert_{LC_{q_{i},\lambda_{i}}^{\left \{  x_{0}\right \}
}}r^{\frac{n}{p}}\int \limits_{2r}^{\infty}\left(  1+\ln \frac{t}{r}\right)
^{2}%
{\displaystyle \prod \limits_{i=1}^{2}}
\Vert f_{i}\Vert_{L_{p_{i}}(B(x_{0},t))}\\
& \times \frac{dt}{t^{n\left(  \left(  \frac{1}{p_{1}}+\frac{1}{p_{2}}\right)
-\left(  \lambda_{1}+\lambda_{2}\right)  \right)  +1}}.
\end{align*}

$\left(  ii\right)  $ For $G_{2}=\left \Vert T_{\left(  b_{1},b_{2}\right)
}^{\left(  2\right)  }\left(  f_{1}^{0},f_{2}^{\infty}\right)  \right \Vert
_{L_{p}\left(  B\left(  x_{0},r\right)  \right)  }$, we also write%
\begin{align*}
G_{2}  & \lesssim \left \Vert \left[  \left(  b_{1}-\left \{  b_{1}\right \}
_{B}\right)  \right]  \left[  \left(  b_{2}-\left \{  b_{2}\right \}
_{B}\right)  \right]  T^{\left(  2\right)  }\left(  f_{1}^{0},f_{2}^{\infty
}\right)  \right \Vert _{L_{p}\left(  B\left(  x_{0},r\right)  \right)  }\\
& +\left \Vert \left[  \left(  b_{1}-\left \{  b_{1}\right \}  _{B}\right)
\right]  T^{\left(  2\right)  }\left[  f_{1}^{0},\left(  b_{2}-\left \{
b_{2}\right \}  _{B}\right)  f_{2}^{\infty}\right]  \right \Vert _{L_{p}\left(
B\left(  x_{0},r\right)  \right)  }\\
& +\left \Vert \left[  \left(  b_{2}-\left \{  b_{2}\right \}  _{B}\right)
\right]  T^{\left(  2\right)  }\left[  \left(  b_{1}-\left \{  b_{1}\right \}
_{B}\right)  f_{1}^{0},f_{2}^{\infty}\right]  \right \Vert _{L_{p}\left(
B\left(  x_{0},r\right)  \right)  }\\
& +\left \Vert T^{\left(  2\right)  }\left[  \left(  b_{1}-\left \{
b_{1}\right \}  _{B}\right)  f_{1}^{0},\left(  b_{2}-\left \{  b_{2}\right \}
_{B}\right)  f_{2}^{\infty}\right]  \right \Vert _{L_{p}\left(  B\left(
x_{0},r\right)  \right)  }\\
& \equiv G_{21}+G_{22}+G_{23}+G_{24}.
\end{align*}

Let $1<\overline{p},\overline{q}<\infty$, such that $\frac{1}{\overline{p}%
}=\frac{1}{p_{1}}+\frac{1}{p_{2}}$ and $\frac{1}{\overline{q}}=\frac{1}{q_{1}%
}+\frac{1}{q_{2}}$. Then, using H\"{o}lder's inequality we have%
\begin{align*}
G_{21}  & =\left \Vert \left[  \left(  b_{1}-\left \{  b_{1}\right \}
_{B}\right)  \right]  \left[  \left(  b_{2}-\left \{  b_{2}\right \}
_{B}\right)  \right]  T^{\left(  2\right)  }\left(  f_{1}^{0},f_{2}^{\infty
}\right)  \right \Vert _{L_{p}\left(  B\left(  x_{0},r\right)  \right)  }\\
& \lesssim \left \Vert \left[  \left(  b_{1}-\left \{  b_{1}\right \}
_{B}\right)  \right]  \left[  \left(  b_{2}-\left \{  b_{2}\right \}
_{B}\right)  \right]  \right \Vert _{L_{\overline{q}}\left(  B\right)
}\left \Vert T^{\left(  2\right)  }\left(  f_{1}^{0},f_{2}^{\infty}\right)
\right \Vert _{L_{\overline{p}}\left(  B\right)  }\\
& \lesssim \left \Vert \left(  b_{1}-\left \{  b_{1}\right \}  _{B}\right)
\right \Vert _{L_{q_{1}}\left(  B\right)  }\left \Vert \left(  b_{2}-\left \{
b_{2}\right \}  _{B}\right)  \right \Vert _{L_{q_{2}}\left(  B\right)  }%
r^{\frac{n}{\overline{p}}}\int \limits_{2r}^{\infty}%
{\displaystyle \prod \limits_{i=1}^{2}}
\Vert f_{i}\Vert_{L_{p_{i}}(B(x_{0},t))}\frac{dt}{t^{\frac{n}{\overline{p}}%
+1}}\\
& \lesssim%
{\displaystyle \prod \limits_{i=1}^{2}}
\Vert \overrightarrow{b}\Vert_{LC_{q_{i},\lambda_{i}}^{\left \{  x_{0}\right \}
}}r^{n\left(  \frac{1}{q_{1}}+\frac{1}{q_{2}}\right)  +n\left(  \lambda
_{1}+\lambda_{2}\right)  }r^{n\left(  \frac{1}{p_{1}}+\frac{1}{p_{2}}\right)
}\int \limits_{2r}^{\infty}\left(  1+\ln \frac{t}{r}\right)  ^{2}%
{\displaystyle \prod \limits_{i=1}^{2}}
\Vert f_{i}\Vert_{L_{p_{i}}(B(x_{0},t))}\frac{dt}{t^{n\left(  \frac{1}{p_{1}%
}+\frac{1}{p_{2}}\right)  +1}}\\
& \lesssim%
{\displaystyle \prod \limits_{i=1}^{2}}
\Vert \overrightarrow{b}\Vert_{LC_{q_{i},\lambda_{i}}^{\left \{  x_{0}\right \}
}}r^{\frac{n}{p}}\int \limits_{2r}^{\infty}\left(  1+\ln \frac{t}{r}\right)
^{2}%
{\displaystyle \prod \limits_{i=1}^{2}}
\Vert f_{i}\Vert_{L_{p_{i}}(B(x_{0},t))}\frac{dt}{t^{n\left(  \left(  \frac
{1}{p_{1}}+\frac{1}{p_{2}}\right)  -\left(  \lambda_{1}+\lambda_{2}\right)
\right)  +1}},
\end{align*}
where in the second inequality we have used the following fact:

It is clear that $\left \vert \left(  x_{0}-y_{1},\text{ }x_{0}-y_{2}\right)
\right \vert ^{2n}\geq \left \vert x_{0}-y_{2}\right \vert ^{2n}$. By the
condition (\ref{4}) with $m=2$ and H\"{o}lder's inequality, we have%
\begin{align*}
\left \vert T^{\left(  2\right)  }\left(  f_{1}^{0},f_{2}^{\infty}\right)
\left(  x\right)  \right \vert  & \lesssim%
{\displaystyle \int \limits_{{\mathbb{R}^{n}}}}
{\displaystyle \int \limits_{{\mathbb{R}^{n}}}}
\frac{\left \vert f_{1}^{0}\left(  y_{1}\right)  \right \vert \left \vert
f_{2}^{\infty}\left(  y_{2}\right)  \right \vert }{\left \vert \left(
x-y_{1},x-y_{2}\right)  \right \vert ^{2n}}dy_{1}dy_{2}\\
& \lesssim \int \limits_{2B}\left \vert f_{1}\left(  y_{1}\right)  \right \vert
dy_{1}\int \limits_{\left(  2B\right)  ^{C}}\frac{\left \vert f_{2}\left(
y_{2}\right)  \right \vert }{\left \vert x_{0}-y_{2}\right \vert ^{2n}}dy_{2}\\
& \approx \int \limits_{2B}\left \vert f_{1}\left(  y_{1}\right)  \right \vert
dy_{1}\int \limits_{\left(  2B\right)  ^{C}}\left \vert f_{2}\left(
y_{2}\right)  \right \vert \int \limits_{\left \vert x_{0}-y_{2}\right \vert
}^{\infty}\frac{dt}{t^{2n+1}}dy_{2}\\
& \lesssim \left \Vert f_{1}\right \Vert _{L_{p_{1}}\left(  2B\right)
}\left \vert 2B\right \vert ^{1-\frac{1}{p_{1}}}%
{\displaystyle \int \limits_{2r}^{\infty}}
\left \Vert f_{2}\right \Vert _{L_{p_{2}}\left(  B\left(  x_{0},t\right)
\right)  }\left \vert B\left(  x_{0},t\right)  \right \vert ^{1-\frac{1}{p_{2}}%
}\frac{dt}{t^{2n+1}}\\
& \lesssim \int \limits_{2r}^{\infty}%
{\displaystyle \prod \limits_{i=1}^{2}}
\Vert f_{i}\Vert_{L_{p_{i}}(B(x_{0},t))}\frac{dt}{t^{\frac{n}{\overline{p}}%
+1}},
\end{align*}
where $\frac{1}{\overline{p}}=\frac{1}{p_{1}}+\frac{1}{p_{2}}$. Thus, the
inequality
\[
\left \Vert T^{\left(  2\right)  }\left(  f_{1}^{0},f_{2}^{\infty}\right)
\right \Vert _{L_{\overline{p}}\left(  B\left(  x_{0},r\right)  \right)
}\lesssim r^{\frac{n}{\overline{p}}}\int \limits_{2r}^{\infty}%
{\displaystyle \prod \limits_{i=1}^{2}}
\Vert f_{i}\Vert_{L_{p_{i}}(B(x_{0},t))}\frac{dt}{t^{\frac{n}{\overline{p}}%
+1}}%
\]
is valid.

On the other hand, for the estimates used in $G_{22}$, $G_{23}$, we have to
prove the below inequality:%
\begin{equation}
\left \vert T^{\left(  2\right)  }\left[  f_{1}^{0},\left(  b_{2}\left(
\cdot \right)  -\left \{  b_{2}\right \}  _{B}\right)  f_{2}^{\infty}\right]
\left(  x\right)  \right \vert \lesssim \Vert b_{2}\Vert_{LC_{q_{2},\lambda_{2}%
}^{\left \{  x_{0}\right \}  }}\int \limits_{2r}^{\infty}\left(  1+\ln \frac{t}%
{r}\right)  ^{2}%
{\displaystyle \prod \limits_{i=1}^{2}}
\Vert f_{i}\Vert_{L_{p_{i}}(B(x_{0},t))}\frac{dt}{t^{n\left(  \left(  \frac
{1}{p_{1}}+\frac{1}{p_{2}}\right)  -\lambda_{2}\right)  +1}}.\label{11*}%
\end{equation}
Indeed, it is clear that $\left \vert \left(  x_{0}-y_{1},\text{ }x_{0}%
-y_{2}\right)  \right \vert ^{2n}\geq \left \vert x_{0}-y_{2}\right \vert ^{2n}$.
Moreover, using the conditions (\ref{10*}) and (\ref{4}) with $m=2$, we have%
\begin{align*}
& \left \vert T^{\left(  2\right)  }\left[  f_{1}^{0},\left(  b_{2}\left(
\cdot \right)  -\left \{  b_{2}\right \}  _{B}\right)  f_{2}^{\infty}\right]
\left(  x\right)  \right \vert \\
& \lesssim%
{\displaystyle \int \limits_{2B}}
\left \vert f_{1}\left(  y_{1}\right)  \right \vert dy_{1}%
{\displaystyle \int \limits_{\left(  2B\right)  ^{C}}}
\frac{\left \vert b_{2}\left(  y_{2}\right)  -\left \{  b_{2}\right \}
_{B}\right \vert \left \vert f_{2}\left(  y_{2}\right)  \right \vert }{\left \vert
x_{0}-y_{2}\right \vert ^{2n}}dy_{2}.
\end{align*}

It's obvious that%
\begin{equation}%
{\displaystyle \int \limits_{2B}}
\left \vert f_{1}\left(  y_{1}\right)  \right \vert dy_{1}\lesssim \left \Vert
f_{1}\right \Vert _{L_{p_{1}}\left(  2B\right)  }\left \vert 2B\right \vert
^{1-\frac{1}{p_{1}}},\label{11}%
\end{equation}
and using H\"{o}lder's inequality and by (\ref{*})%
\begin{align*}
&
{\displaystyle \int \limits_{\left(  2B\right)  ^{C}}}
\frac{\left \vert b_{2}\left(  y_{2}\right)  -\left \{  b_{2}\right \}
_{B}\right \vert \left \vert f_{2}\left(  y_{2}\right)  \right \vert }{\left \vert
x_{0}-y_{2}\right \vert ^{2n}}dy_{2}\\
& \lesssim%
{\displaystyle \int \limits_{\left(  2B\right)  ^{C}}}
\left \vert b_{2}\left(  y_{2}\right)  -\left \{  b_{2}\right \}  _{B}\right \vert
\left \vert f_{2}\left(  y_{2}\right)  \right \vert \left[
{\displaystyle \int \limits_{\left \vert x_{0}-y_{2}\right \vert }^{\infty}}
\frac{dt}{t^{2n+1}}\right]  dy_{2}\\
& \lesssim%
{\displaystyle \int \limits_{2r}^{\infty}}
\left \Vert b_{2}\left(  y_{2}\right)  -\left \{  b_{2}\right \}  _{B\left(
x_{0},t\right)  }\right \Vert _{L_{q_{2}}\left(  B\left(  x_{0},t\right)
\right)  }\left \Vert f_{2}\right \Vert _{L_{p_{2}}\left(  B\left(
x_{0},t\right)  \right)  }\left \vert B\left(  x_{0},t\right)  \right \vert
^{1-\left(  \frac{1}{p_{1}}+\frac{1}{p_{2}}\right)  }\frac{dt}{t^{2n+1}}\\
& +\left \vert \left \{  b_{2}\right \}  _{B\left(  x_{0},t\right)  }-\left \{
b_{2}\right \}  _{B\left(  x_{0},r\right)  }\right \vert \left \Vert
f_{2}\right \Vert _{L_{p_{2}}\left(  B\left(  x_{0},t\right)  \right)
}\left \vert B\left(  x_{0},t\right)  \right \vert ^{1-\frac{1}{p_{2}}}\frac
{dt}{t^{2n+1}}\\
& \lesssim \Vert b_{2}\Vert_{LC_{q_{2},\lambda_{2}}^{\left \{  x_{0}\right \}  }}%
{\displaystyle \int \limits_{2r}^{\infty}}
\left \vert B\left(  x_{0},t\right)  \right \vert ^{\frac{1}{q_{2}}+\lambda_{2}%
}\left \Vert f_{2}\right \Vert _{L_{p_{2}}\left(  B\left(  x_{0},t\right)
\right)  }\left \vert B\left(  x_{0},t\right)  \right \vert ^{1-\left(  \frac
{1}{p_{2}}+\frac{1}{q_{2}}\right)  }\frac{dt}{t^{2n+1}}\\
& \lesssim \Vert b_{2}\Vert_{LC_{q_{2},\lambda_{2}}^{\left \{  x_{0}\right \}  }%
}\int \limits_{2r}^{\infty}\left(  1+\ln \frac{t}{r}\right)  \left \vert B\left(
x_{0},t\right)  \right \vert ^{\lambda_{2}}\left \Vert f_{2}\right \Vert
_{L_{p_{2}}\left(  B\left(  x_{0},t\right)  \right)  }\left \vert B\left(
x_{0},t\right)  \right \vert ^{1-\frac{1}{p_{2}}}\frac{dt}{t^{2n+1}}%
\end{align*}%
\begin{equation}
\lesssim \Vert b_{2}\Vert_{LC_{q_{2},\lambda_{2}}^{\left \{  x_{0}\right \}  }%
}\int \limits_{2r}^{\infty}\left(  1+\ln \frac{t}{r}\right)  ^{2}\left \Vert
f_{2}\right \Vert _{L_{p_{2}}\left(  B\left(  x_{0},t\right)  \right)  }%
\frac{dt}{t^{n\left(  1+\frac{1}{p_{2}}-\lambda_{2}\right)  +1}}.\label{F}%
\end{equation}
Hence, by (\ref{11}) and (\ref{F}), it follows that:%
\begin{align*}
& \left \vert T^{\left(  2\right)  }\left[  f_{1}^{0},\left(  b_{2}\left(
\cdot \right)  -\left \{  b_{2}\right \}  _{B}\right)  f_{2}^{\infty}\right]
\left(  x\right)  \right \vert \\
& \lesssim \Vert b_{2}\Vert_{LC_{q_{2},\lambda_{2}}^{\left \{  x_{0}\right \}  }%
}\left \Vert f_{1}\right \Vert _{L_{p_{1}}\left(  2B\right)  }\left \vert
2B\right \vert ^{1-\frac{1}{p_{1}}}\int \limits_{2r}^{\infty}\left(  1+\ln
\frac{t}{r}\right)  ^{2}\left \Vert f_{2}\right \Vert _{L_{p_{2}}\left(
B\left(  x_{0},t\right)  \right)  }\frac{dt}{t^{n\left(  1+\frac{1}{p_{2}%
}-\lambda_{2}\right)  +1}}\\
& \lesssim \Vert b_{2}\Vert_{LC_{q_{2},\lambda_{2}}^{\left \{  x_{0}\right \}  }%
}\int \limits_{2r}^{\infty}\left(  1+\ln \frac{t}{r}\right)  ^{2}%
{\displaystyle \prod \limits_{i=1}^{2}}
\Vert f_{i}\Vert_{L_{p_{i}}(B(x_{0},t))}\frac{dt}{t^{n\left(  \left(  \frac
{1}{p_{1}}+\frac{1}{p_{2}}\right)  -\lambda_{2}\right)  +1}}.
\end{align*}
This completes the proof of inequality (\ref{11*}).

Thus, let $1<\tau<\infty$, such that $\frac{1}{p}=\frac{1}{q_{1}}+\frac
{1}{\tau}$. Then, to estimate $G_{22}$, similar to the estimates for $G_{11}$,
using H\"{o}lder's inequality and from (\ref{F}), we get%
\begin{align*}
G_{22}  & =\left \Vert \left[  \left(  b_{1}-\left \{  b_{1}\right \}
_{B}\right)  \right]  T^{\left(  2\right)  }\left[  f_{1}^{0},\left(
b_{2}-\left \{  b_{2}\right \}  _{B}\right)  f_{2}^{\infty}\right]  \right \Vert
_{L_{p}\left(  B\left(  x_{0},r\right)  \right)  }\\
& \lesssim \left \Vert \left(  b_{1}-\left \{  b_{1}\right \}  _{B}\right)
\right \Vert _{L_{q_{1}}\left(  B\right)  }\left \Vert T^{\left(  2\right)  }
\left[  f_{1}^{0},\left(  b_{2}-\left \{  b_{2}\right \}  _{B}\right)
f_{2}^{\infty}\right]  \right \Vert _{L_{\tau}\left(  B\right)  }\\
& \lesssim%
{\displaystyle \prod \limits_{i=1}^{2}}
\Vert \overrightarrow{b}\Vert_{LC_{q_{i},\lambda_{i}}^{\left \{  x_{0}\right \}
}}\left \vert B\right \vert ^{\lambda_{1}+\frac{1}{q_{1}}+\frac{1}{\tau}}%
\int \limits_{2r}^{\infty}\left(  1+\ln \frac{t}{r}\right)  ^{2}%
{\displaystyle \prod \limits_{i=1}^{2}}
\Vert f_{i}\Vert_{L_{p_{i}}(B(x_{0},t))}\frac{dt}{t^{n\left(  \left(  \frac
{1}{p_{1}}+\frac{1}{p_{2}}\right)  -\lambda_{2}\right)  +1}}\\
& \lesssim%
{\displaystyle \prod \limits_{i=1}^{2}}
\Vert \overrightarrow{b}\Vert_{LC_{q_{i},\lambda_{i}}^{\left \{  x_{0}\right \}
}}r^{\frac{n}{p}}\int \limits_{2r}^{\infty}\left(  1+\ln \frac{t}{r}\right)
^{2}%
{\displaystyle \prod \limits_{i=1}^{2}}
\Vert f_{i}\Vert_{L_{p_{i}}(B(x_{0},t))}\frac{dt}{t^{n\left(  \left(  \frac
{1}{p_{1}}+\frac{1}{p_{2}}\right)  -\left(  \lambda_{1}+\lambda_{2}\right)
\right)  +1}}.
\end{align*}
Similarly, $G_{23}$ has the same estimate above, here we omit the details,
thus the inequality%
\begin{align*}
G_{23}  & =\left \Vert \left[  \left(  b_{2}-\left \{  b_{2}\right \}
_{B}\right)  \right]  T^{\left(  2\right)  }\left[  \left(  b_{1}-\left \{
b_{1}\right \}  _{B}\right)  f_{1}^{0},f_{2}^{\infty}\right]  \right \Vert
_{L_{p}\left(  B\left(  x_{0},r\right)  \right)  }\\
& \lesssim%
{\displaystyle \prod \limits_{i=1}^{2}}
\Vert \overrightarrow{b}\Vert_{LC_{q_{i},\lambda_{i}}^{\left \{  x_{0}\right \}
}}r^{\frac{n}{p}}\int \limits_{2r}^{\infty}\left(  1+\ln \frac{t}{r}\right)
^{2}%
{\displaystyle \prod \limits_{i=1}^{2}}
\Vert f_{i}\Vert_{L_{p_{i}}(B(x_{0},t))}\frac{dt}{t^{n\left(  \left(  \frac
{1}{p_{1}}+\frac{1}{p_{2}}\right)  -\left(  \lambda_{1}+\lambda_{2}\right)
\right)  +1}}%
\end{align*}
is valid.

Now we turn to estimate $G_{24}$. Similar to (\ref{F}), we have to prove the
following estimate for $G_{24}$:%
\begin{align}
& \left \vert T^{\left(  2\right)  }\left[  \left(  b_{1}-\left \{
b_{1}\right \}  _{B}\right)  f_{1}^{0},\left(  b_{2}-\left \{  b_{2}\right \}
_{B}\right)  f_{2}^{\infty}\right]  \left(  x\right)  \right \vert \nonumber \\
& \lesssim%
{\displaystyle \prod \limits_{i=1}^{2}}
\Vert \overrightarrow{b}\Vert_{LC_{q_{i},\lambda_{i}}^{\left \{  x_{0}\right \}
}}\int \limits_{2r}^{\infty}\left(  1+\ln \frac{t}{r}\right)  ^{2}%
{\displaystyle \prod \limits_{i=1}^{2}}
\Vert f_{i}\Vert_{L_{p_{i}}(B(x_{0},t))}\frac{dt}{t^{n\left(  \left(  \frac
{1}{p_{1}}+\frac{1}{p_{2}}\right)  -\left(  \lambda_{1}+\lambda_{2}\right)
\right)  +1}}.\label{G}%
\end{align}
Firstly, using the condition (\ref{4}) with $m=2$, we have%
\begin{align*}
& \left \vert T^{\left(  2\right)  }\left[  \left(  b_{1}-\left \{
b_{1}\right \}  _{B}\right)  f_{1}^{0},\left(  b_{2}-\left \{  b_{2}\right \}
_{B}\right)  f_{2}^{\infty}\right]  \left(  x\right)  \right \vert \\
& \lesssim%
{\displaystyle \int \limits_{2B}}
\left \vert b_{1}\left(  y_{1}\right)  -\left \{  b_{1}\right \}  _{B}\right \vert
\left \vert f_{1}\left(  y_{1}\right)  \right \vert dy_{1}%
{\displaystyle \int \limits_{\left(  2B\right)  ^{C}}}
\frac{\left \vert b_{2}\left(  y_{2}\right)  -\left \{  b_{2}\right \}
_{B}\right \vert \left \vert f_{2}\left(  y_{2}\right)  \right \vert }{\left \vert
x_{0}-y_{2}\right \vert ^{2n}}dy_{2}.
\end{align*}

It's obvious that%
\begin{equation}%
{\displaystyle \int \limits_{2B}}
\left \vert b_{1}\left(  y_{1}\right)  -\left \{  b_{1}\right \}  _{B}\right \vert
\left \vert f_{1}\left(  y_{1}\right)  \right \vert dy_{1}\lesssim \Vert
b_{1}\Vert_{LC_{q_{1},\lambda_{1}}^{\left \{  x_{0}\right \}  }}\left \vert
B\right \vert ^{\lambda_{1}+1-\frac{1}{p_{1}}}\left \Vert f_{1}\right \Vert
_{L_{p_{1}}\left(  2B\right)  }.\label{12}%
\end{equation}
Then, by (\ref{F}) and (\ref{12}) we get (\ref{G}). This completes the proof
of inequality (\ref{G}). Therefore, by (\ref{G}) we deduce that%
\begin{align*}
G_{24}  & =\left \Vert T^{\left(  2\right)  }\left[  \left(  b_{1}-\left \{
b_{1}\right \}  _{B}\right)  f_{1}^{0},\left(  b_{2}-\left \{  b_{2}\right \}
_{B}\right)  f_{2}^{\infty}\right]  \right \Vert _{L_{p}\left(  B\left(
x_{0},r\right)  \right)  }\\
& \lesssim%
{\displaystyle \prod \limits_{i=1}^{2}}
\Vert \overrightarrow{b}\Vert_{LC_{q_{i},\lambda_{i}}^{\left \{  x_{0}\right \}
}}r^{\frac{n}{p}}\int \limits_{2r}^{\infty}\left(  1+\ln \frac{t}{r}\right)
^{2}%
{\displaystyle \prod \limits_{i=1}^{2}}
\Vert f_{i}\Vert_{L_{p_{i}}(B(x_{0},t))}\frac{dt}{t^{n\left(  \left(  \frac
{1}{p_{1}}+\frac{1}{p_{2}}\right)  -\left(  \lambda_{1}+\lambda_{2}\right)
\right)  +1}}.
\end{align*}
Considering estimates $G_{21},$ $G_{22}$, $G_{23}$, $G_{24}$ together, we get
the desired conclusion%
\begin{align*}
G_{2}  & =\left \Vert T_{\left(  b_{1},b_{2}\right)  }^{\left(  2\right)
}\left(  f_{1}^{0},f_{2}^{\infty}\right)  \right \Vert _{L_{p}\left(  B\left(
x_{0},r\right)  \right)  }\lesssim%
{\displaystyle \prod \limits_{i=1}^{2}}
\Vert \overrightarrow{b}\Vert_{LC_{q_{i},\lambda_{i}}^{\left \{  x_{0}\right \}
}}r^{\frac{n}{p}}\int \limits_{2r}^{\infty}\left(  1+\ln \frac{t}{r}\right)
^{2}%
{\displaystyle \prod \limits_{i=1}^{2}}
\Vert f_{i}\Vert_{L_{p_{i}}(B(x_{0},t))}\\
& \times \frac{dt}{t^{n\left(  \left(  \frac{1}{p_{1}}+\frac{1}{p_{2}}\right)
-\left(  \lambda_{1}+\lambda_{2}\right)  \right)  +1}}.
\end{align*}

Similar to $G_{2}$, we can also get the estimates for $F_{3}$,
\begin{align*}
G_{3}  & =\left \Vert T_{\left(  b_{1},b_{2}\right)  }^{\left(  2\right)
}\left(  f_{1}^{\infty},f_{2}^{0}\right)  \right \Vert _{L_{p}\left(  B\left(
x_{0},r\right)  \right)  }\lesssim%
{\displaystyle \prod \limits_{i=1}^{2}}
\Vert \overrightarrow{b}\Vert_{LC_{q_{i},\lambda_{i}}^{\left \{  x_{0}\right \}
}}r^{\frac{n}{p}}\int \limits_{2r}^{\infty}\left(  1+\ln \frac{t}{r}\right)
^{2}%
{\displaystyle \prod \limits_{i=1}^{2}}
\Vert f_{i}\Vert_{L_{p_{i}}(B(x_{0},t))}\\
& \times \frac{dt}{t^{n\left(  \left(  \frac{1}{p_{1}}+\frac{1}{p_{2}}\right)
-\left(  \lambda_{1}+\lambda_{2}\right)  \right)  +1}}.
\end{align*}
Finally, for $G_{4}=\left \Vert T_{\left(  b_{1},b_{2}\right)  }^{\left(
2\right)  }\left(  f_{1}^{\infty},f_{2}^{\infty}\right)  \right \Vert
_{L_{p}\left(  B\left(  x_{0},r\right)  \right)  }$, we write%
\begin{align*}
G_{4}  & \lesssim \left \Vert \left[  \left(  b_{1}-\left \{  b_{1}\right \}
_{B}\right)  \right]  \left[  \left(  b_{2}-\left \{  b_{2}\right \}
_{B}\right)  \right]  T^{\left(  2\right)  }\left(  f_{1}^{\infty}%
,f_{2}^{\infty}\right)  \right \Vert _{L_{p}\left(  B\left(  x_{0},r\right)
\right)  }\\
& +\left \Vert \left[  \left(  b_{1}-\left \{  b_{1}\right \}  _{B}\right)
\right]  T^{\left(  2\right)  }\left[  f_{1}^{\infty},\left(  b_{2}-\left \{
b_{2}\right \}  _{B}\right)  f_{2}^{\infty}\right]  \right \Vert _{L_{p}\left(
B\left(  x_{0},r\right)  \right)  }\\
& +\left \Vert \left[  \left(  b_{2}-\left \{  b_{2}\right \}  _{B}\right)
\right]  T^{\left(  2\right)  }\left[  \left(  b_{1}-\left \{  b_{1}\right \}
_{B}\right)  f_{1}^{\infty},f_{2}^{\infty}\right]  \right \Vert _{L_{p}\left(
B\left(  x_{0},r\right)  \right)  }\\
& +\left \Vert T^{\left(  2\right)  }\left[  \left(  b_{1}-\left \{
b_{1}\right \}  _{B}\right)  f_{1}^{\infty},\left(  b_{2}-\left \{
b_{2}\right \}  _{B}\right)  f_{2}^{\infty}\right]  \right \Vert _{L_{p}\left(
B\left(  x_{0},r\right)  \right)  }\\
& \equiv G_{41}+G_{42}+G_{43}+G_{44}.
\end{align*}

Now, let us estimate $G_{41}$, $G_{42}$, $G_{43}$, $G_{44}$, respectively.

For the term $G_{41}$, let $1<\tau<\infty$, such that $\frac{1}{p}=\left(
\frac{1}{q_{1}}+\frac{1}{q_{2}}\right)  +\frac{1}{\tau}$, $\frac{1}{\tau
}=\frac{1}{p_{1}}+\frac{1}{p_{2}}$. Then, by H\"{o}lder's inequality we get%
\begin{align*}
G_{41}  & =\left \Vert \left[  \left(  b_{1}-\left \{  b_{1}\right \}
_{B}\right)  \right]  \left[  \left(  b_{2}-\left \{  b_{2}\right \}
_{B}\right)  \right]  T^{\left(  2\right)  }\left(  f_{1}^{\infty}%
,f_{2}^{\infty}\right)  \right \Vert _{L_{p}\left(  B\left(  x_{0},r\right)
\right)  }\\
& \lesssim \left \Vert \left(  b_{1}-\left \{  b_{1}\right \}  _{B}\right)
\right \Vert _{L_{q_{1}}\left(  B\right)  }\left \Vert \left(  b_{2}-\left \{
b_{2}\right \}  _{B}\right)  \right \Vert _{L_{q_{2}}\left(  B\right)
}\left \Vert T^{\left(  2\right)  }\left(  f_{1}^{\infty},f_{2}^{\infty
}\right)  \right \Vert _{L_{\tau}\left(  B\right)  }\\
& \lesssim%
{\displaystyle \prod \limits_{i=1}^{2}}
\Vert \overrightarrow{b}\Vert_{LC_{q_{i},\lambda_{i}}^{\left \{  x_{0}\right \}
}}\left \vert B\right \vert ^{\left(  \lambda_{1}+\lambda_{2}\right)  +\left(
\frac{1}{q_{1}}+\frac{1}{q_{2}}\right)  }r^{\frac{n}{\tau}}\int \limits_{2r}%
^{\infty}\left(  1+\ln \frac{t}{r}\right)  ^{2}%
{\displaystyle \prod \limits_{i=1}^{2}}
\Vert f_{i}\Vert_{L_{p_{i}}(B(x_{0},t))}\frac{dt}{t^{\frac{n}{\tau}+1}}\\
& \lesssim%
{\displaystyle \prod \limits_{i=1}^{2}}
\Vert \overrightarrow{b}\Vert_{LC_{q_{i},\lambda_{i}}^{\left \{  x_{0}\right \}
}}r^{\frac{n}{p}}\int \limits_{2r}^{\infty}\left(  1+\ln \frac{t}{r}\right)
^{2}%
{\displaystyle \prod \limits_{i=1}^{2}}
\Vert f_{i}\Vert_{L_{p_{i}}(B(x_{0},t))}\frac{dt}{t^{n\left(  \left(  \frac
{1}{p_{1}}+\frac{1}{p_{2}}\right)  -\left(  \lambda_{1}+\lambda_{2}\right)
\right)  +1}},
\end{align*}
where in the second inequality we have used the following fact:

Noting that $\left \vert \left(  x_{0}-y_{1},\text{ }x_{0}-y_{2}\right)
\right \vert ^{2n}\geq \left \vert x_{0}-y_{1}\right \vert ^{n}\left \vert
x_{0}-y_{2}\right \vert ^{n}$. Using the condition (\ref{4}) with $m=2$ and by
H\"{o}lder's inequality, we get%
\begin{align*}
& \left \vert T_{\alpha}^{\left(  2\right)  }\left(  f_{1}^{\infty}%
,f_{2}^{\infty}\right)  \left(  x\right)  \right \vert \\
& \lesssim%
{\displaystyle \int \limits_{\mathbb{R} ^{n}}}
{\displaystyle \int \limits_{\mathbb{R} ^{n}}}
\frac{\left \vert f_{1}\left(  y_{1}\right)  \chi_{\left(  2B\right)  ^{c}%
}\right \vert \left \vert f_{2}\left(  y_{2}\right)  \chi_{\left(  2B\right)
^{c}}\right \vert }{\left \vert \left(  x_{0}-y_{1},x_{0}-y_{2}\right)
\right \vert ^{2n}}dy_{1}dy_{2}\\
& \lesssim%
{\displaystyle \int \limits_{\left(  2B\right)  ^{c}}}
{\displaystyle \int \limits_{\left(  2B\right)  ^{c}}}
\frac{\left \vert f_{1}\left(  y_{1}\right)  \right \vert \left \vert
f_{2}\left(  y_{2}\right)  \right \vert }{\left \vert x_{0}-y_{1}\right \vert
^{n}\left \vert x_{0}-y_{2}\right \vert ^{n}}dy_{1}dy_{2}\\
& \lesssim%
{\displaystyle \sum \limits_{j=1}^{\infty}}
{\displaystyle \prod \limits_{i=1}^{2}}
{\displaystyle \int \limits_{B\left(  x_{0},2^{j+1}r\right)  \backslash B\left(
x_{0},2^{j}r\right)  }}
\frac{\left \vert f_{i}\left(  y_{i}\right)  \right \vert }{\left \vert
x_{0}-y_{i}\right \vert ^{n}}dy_{i}\\
& \lesssim%
{\displaystyle \sum \limits_{j=1}^{\infty}}
{\displaystyle \prod \limits_{i=1}^{2}}
\left(  2^{j}r\right)  ^{-n}%
{\displaystyle \int \limits_{B\left(  x_{0},2^{j+1}r\right)  }}
\left \vert f_{i}\left(  y_{i}\right)  \right \vert dy_{i}\\
& \lesssim%
{\displaystyle \sum \limits_{j=1}^{\infty}}
\left(  2^{j}r\right)  ^{-2n}%
{\displaystyle \prod \limits_{i=1}^{2}}
\left \Vert f_{i}\right \Vert _{L_{p_{i}}(B\left(  x_{0},2^{j+1}r\right)
)}\left \vert B\left(  x_{0},2^{j+1}r\right)  \right \vert ^{1-\frac{1}{p_{i}}%
}\\
& \lesssim%
{\displaystyle \sum \limits_{j=1}^{\infty}}
{\displaystyle \int \limits_{2^{j+1}r}^{2^{j+2}r}}
\left(  2^{j+1}r\right)  ^{-2n-1}%
{\displaystyle \prod \limits_{i=1}^{2}}
\left \Vert f_{i}\right \Vert _{L_{p_{i}}(B\left(  x_{0},2^{j+1}r\right)
)}\left \vert B\left(  x_{0},2^{j+1}r\right)  \right \vert ^{1-\frac{1}{p_{i}}%
}dt\\
& \lesssim%
{\displaystyle \sum \limits_{j=1}^{\infty}}
{\displaystyle \int \limits_{2^{j+1}r}^{2^{j+2}r}}
{\displaystyle \prod \limits_{i=1}^{2}}
\left \Vert f_{i}\right \Vert _{L_{p_{i}}(B\left(  x_{0},t\right)  )}\left \vert
B\left(  x_{0},t\right)  \right \vert ^{1-\frac{1}{p_{i}}}\frac{dt}{t^{2n+1}}\\
& \lesssim%
{\displaystyle \int \limits_{2r}^{\infty}}
\prod \limits_{i=1}^{2}\left \Vert f_{i}\right \Vert _{L_{p_{i}}\left(  B\left(
x_{0},t\right)  \right)  }\left \vert B\left(  x_{0},t\right)  \right \vert
^{2-\left(  \frac{1}{p_{1}}+\frac{1}{p_{2}}\right)  }\frac{dt}{t^{2n+1}}\\
& \lesssim%
{\displaystyle \int \limits_{2r}^{\infty}}
\left \Vert f_{1}\right \Vert _{L_{p_{1}}\left(  B\left(  x_{0},t\right)
\right)  }\left \Vert f_{2}\right \Vert _{L_{p_{2}}\left(  B\left(
x_{0},t\right)  \right)  }\frac{dt}{t^{\frac{n}{\tau}+1}}.
\end{align*}
where $\frac{1}{\tau}=\frac{1}{p_{1}}+\frac{1}{p_{2}}$. Thus, for $p_{1}$,
$p_{2}\in \left[  1,\infty \right)  $ the inequality
\[
\left \Vert T^{\left(  2\right)  }\left(  f_{1}^{\infty},f_{2}^{\infty}\right)
\right \Vert _{L_{p}\left(  B\left(  x_{0},r\right)  \right)  }\lesssim
r^{\frac{n}{\tau}}\int \limits_{2r}^{\infty}%
{\displaystyle \prod \limits_{i=1}^{2}}
\Vert f_{i}\Vert_{L_{p_{i}}(B(x_{0},t))}\frac{dt}{t^{\frac{n}{\tau}+1}}%
\]
is valid.

For the terms $G_{42}$, $G_{43}$, similar to the estimates used for
(\ref{11*}), we have to prove the following inequality:%
\begin{equation}
\left \vert T^{\left(  2\right)  }\left[  f_{1}^{\infty},\left(  b_{2}-\left \{
b_{2}\right \}  _{B}\right)  f_{2}^{\infty}\right]  \left(  x\right)
\right \vert \lesssim \Vert b_{2}\Vert_{LC_{q_{2},\lambda_{2}}^{\left \{
x_{0}\right \}  }}\int \limits_{2r}^{\infty}\left(  1+\ln \frac{t}{r}\right)
^{2}%
{\displaystyle \prod \limits_{i=1}^{2}}
\Vert f_{i}\Vert_{L_{p_{i}}(B(x_{0},t))}\frac{dt}{t^{n\left(  \left(  \frac
{1}{p_{1}}+\frac{1}{p_{2}}\right)  -\lambda_{2}\right)  +1}}.\label{f3}%
\end{equation}

Indeed, noting that $\left \vert \left(  x_{0}-y_{1},\text{ }x_{0}%
-y_{2}\right)  \right \vert ^{2n}\geq \left \vert x_{0}-y_{1}\right \vert
^{n}\left \vert x_{0}-y_{2}\right \vert ^{n}$. Recalling the estimates used for
$G_{22}$, $G_{23}$, $G_{24}$ and also using the condition (\ref{4}) with
$m=2$, we have%
\begin{align*}
& \left \vert T^{\left(  2\right)  }\left[  f_{1}^{\infty},\left(
b_{2}-\left \{  b_{2}\right \}  _{B}\right)  f_{2}^{\infty}\right]  \left(
x\right)  \right \vert \\
& \lesssim%
{\displaystyle \int \limits_{\mathbb{R} ^{n}}}
{\displaystyle \int \limits_{\mathbb{R} ^{n}}}
\frac{\left \vert b_{2}\left(  y_{2}\right)  -\left \{  b_{2}\right \}
_{B}\right \vert \left \vert f_{1}\left(  y_{1}\right)  \chi_{\left(  2B\right)
^{C}}\right \vert \left \vert f_{2}\left(  y_{2}\right)  \chi_{\left(
2B\right)  ^{C}}\right \vert }{\left \vert \left(  x_{0}-y_{1},\text{ }%
x_{0}-y_{2}\right)  \right \vert ^{2n}}dy_{1}dy_{2}\\
& \lesssim%
{\displaystyle \int \limits_{\left(  2B\right)  ^{C}}}
{\displaystyle \int \limits_{\left(  2B\right)  ^{C}}}
\frac{\left \vert b_{2}\left(  y_{2}\right)  -\left \{  b_{2}\right \}
_{B}\right \vert \left \vert f_{1}\left(  y_{1}\right)  \right \vert \left \vert
f_{2}\left(  y_{2}\right)  \right \vert }{\left \vert x_{0}-y_{1}\right \vert
^{n}\left \vert x_{0}-y_{2}\right \vert ^{n}}dy_{1}dy_{2}\\
& \lesssim%
{\displaystyle \sum \limits_{j=1}^{\infty}}
{\displaystyle \int \limits_{B\left(  x_{0},2^{j+1}r\right)  \backslash B\left(
x_{0},2^{j}r\right)  }}
\frac{\left \vert f_{1}\left(  y_{1}\right)  \right \vert }{\left \vert
x_{0}-y_{1}\right \vert ^{n}}dy_{1}%
{\displaystyle \int \limits_{B\left(  x_{0},2^{j+1}r\right)  \backslash B\left(
x_{0},2^{j}r\right)  }}
\frac{\left \vert b_{2}\left(  y_{2}\right)  -\left \{  b_{2}\right \}
_{B}\right \vert \left \vert f_{2}\left(  y_{2}\right)  \right \vert }{\left \vert
x_{0}-y_{2}\right \vert ^{n}}dy_{2}\\
& \lesssim%
{\displaystyle \sum \limits_{j=1}^{\infty}}
\left(  2^{j}r\right)  ^{-2n}%
{\displaystyle \int \limits_{B\left(  x_{0},2^{j+1}r\right)  }}
\left \vert f_{1}\left(  y_{1}\right)  \right \vert dy_{1}%
{\displaystyle \int \limits_{B\left(  x_{0},2^{j+1}r\right)  }}
\left \vert b_{2}\left(  y_{2}\right)  -\left \{  b_{2}\right \}  _{B}\right \vert
\left \vert f_{2}\left(  y_{2}\right)  \right \vert dy_{2}.
\end{align*}
It's obvious that%
\begin{equation}%
{\displaystyle \int \limits_{B\left(  x_{0},2^{j+1}r\right)  }}
\left \vert f_{1}\left(  y_{1}\right)  \right \vert dy_{1}\lesssim \left \Vert
f_{1}\right \Vert _{L_{p_{1}}\left(  B\left(  x_{0},2^{j+1}r\right)  \right)
}\left \vert B\left(  x_{0},2^{j+1}r\right)  \right \vert ^{1-\frac{1}{p_{1}}%
},\label{f1}%
\end{equation}
and using H\"{o}lder's inequality and by (\ref{*})
\begin{align*}
&
{\displaystyle \int \limits_{B\left(  x_{0},2^{j+1}r\right)  }}
\left \vert b_{2}\left(  y_{2}\right)  -\left \{  b_{2}\right \}  _{B}\right \vert
\left \vert f_{2}\left(  y_{2}\right)  \right \vert dy_{2}\\
& \lesssim \left \Vert b_{2}\left(  y_{2}\right)  -\left \{  b_{2}\right \}
_{B\left(  x_{0},2^{j+1}r\right)  }\right \Vert _{L_{q_{2}}\left(  B\left(
x_{0},2^{j+1}r\right)  \right)  }\left \Vert f_{2}\right \Vert _{L_{p_{2}%
}\left(  B\left(  x_{0},2^{j+1}r\right)  \right)  }\left \vert B\left(
x_{0},2^{j+1}r\right)  \right \vert ^{1-\left(  \frac{1}{p_{2}}+\frac{1}{q_{2}%
}\right)  }\\
& +\left \vert \left \{  b_{2}\right \}  _{B\left(  x_{0},2^{j+1}r\right)
}-\left \{  b_{2}\right \}  _{B\left(  x_{0},r\right)  }\right \vert \left \Vert
f_{2}\right \Vert _{L_{p_{2}}\left(  B\left(  x_{0},2^{j+1}r\right)  \right)
}\left \vert B\left(  x_{0},2^{j+1}r\right)  \right \vert ^{1-\frac{1}{p_{2}}}\\
& \lesssim \Vert b_{2}\Vert_{LC_{q_{2},\lambda_{2}}^{\left \{  x_{0}\right \}  }%
}\left \vert B\left(  x_{0},2^{j+1}r\right)  \right \vert ^{\frac{1}{q_{2}%
}+\lambda_{2}}\left \Vert f_{2}\right \Vert _{L_{p_{2}}\left(  B\left(
x_{0},2^{j+1}r\right)  \right)  }\left \vert B\left(  x_{0},2^{j+1}r\right)
\right \vert ^{1-\left(  \frac{1}{p_{2}}+\frac{1}{q_{2}}\right)  }\\
& +\Vert b_{2}\Vert_{LC_{q_{2},\lambda_{2}}^{\left \{  x_{0}\right \}  }}\left(
1+\ln \frac{2^{j+1}r}{r}\right)  \left \vert B\left(  x_{0},2^{j+1}r\right)
\right \vert ^{\lambda_{2}}\left \Vert f_{2}\right \Vert _{L_{p_{2}}\left(
B\left(  x_{0},2^{j+1}r\right)  \right)  }\left \vert B\left(  x_{0}%
,2^{j+1}r\right)  \right \vert ^{1-\frac{1}{p_{2}}}%
\end{align*}%
\begin{equation}
\lesssim \Vert b_{2}\Vert_{LC_{q_{2},\lambda_{2}}^{\left \{  x_{0}\right \}  }%
}\int \limits_{2r}^{\infty}\left(  1+\ln \frac{2^{j+1}r}{r}\right)
^{2}\left \vert B\left(  x_{0},2^{j+1}r\right)  \right \vert ^{\lambda_{2}%
-\frac{1}{p_{2}}+1}\left \Vert f_{2}\right \Vert _{L_{p_{2}}\left(  B\left(
x_{0},2^{j+1}r\right)  \right)  }.\label{f2}%
\end{equation}

Hence, by (\ref{f1}) and (\ref{f2}), it follows that:%
\begin{align*}
& \left \vert T^{\left(  2\right)  }\left[  f_{1}^{\infty},\left(
b_{2}-\left \{  b_{2}\right \}  _{B}\right)  f_{2}^{\infty}\right]  \left(
x\right)  \right \vert \\
& \lesssim%
{\displaystyle \sum \limits_{j=1}^{\infty}}
\left(  2^{j}r\right)  ^{-2n}%
{\displaystyle \int \limits_{B\left(  x_{0},2^{j+1}r\right)  }}
\left \vert f_{1}\left(  y_{1}\right)  \right \vert dy_{1}%
{\displaystyle \int \limits_{B\left(  x_{0},2^{j+1}r\right)  }}
\left \vert b_{2}\left(  y_{2}\right)  -\left \{  b_{2}\right \}  _{B}\right \vert
\left \vert f_{2}\left(  y_{2}\right)  \right \vert dy_{2}\\
& \lesssim \Vert b_{2}\Vert_{LC_{q_{2},\lambda_{2}}^{\left \{  x_{0}\right \}  }}%
{\displaystyle \sum \limits_{j=1}^{\infty}}
\left(  2^{j}r\right)  ^{-2n}\left(  1+\ln \frac{2^{j+1}r}{r}\right)
^{2}\left \vert B\left(  x_{0},2^{j+1}r\right)  \right \vert ^{\lambda
_{2}-\left(  \frac{1}{p_{1}}+\frac{1}{p_{2}}\right)  +2}\prod \limits_{i=1}%
^{2}\left \Vert f_{i}\right \Vert _{L_{p_{i}}\left(  B\left(  x_{0}%
,2^{j+1}r\right)  \right)  }\\
& \lesssim \Vert b_{2}\Vert_{LC_{q_{2},\lambda_{2}}^{\left \{  x_{0}\right \}  }}%
{\displaystyle \sum \limits_{j=1}^{\infty}}
\int \limits_{2^{j+1}r}^{2^{j+2}r}\left(  2^{j+1}r\right)  ^{-2n-1}\left(
1+\ln \frac{2^{j+1}r}{r}\right)  ^{2}\left \vert B\left(  x_{0},2^{j+1}r\right)
\right \vert ^{\lambda_{2}-\left(  \frac{1}{p_{1}}+\frac{1}{p_{2}}\right)
+2}\\
& \times \prod \limits_{i=1}^{2}\left \Vert f_{i}\right \Vert _{L_{p_{i}}\left(
B\left(  x_{0},2^{j+1}r\right)  \right)  }dt\\
& \lesssim \Vert b_{2}\Vert_{LC_{q_{2},\lambda_{2}}^{\left \{  x_{0}\right \}  }}%
{\displaystyle \sum \limits_{j=1}^{\infty}}
\int \limits_{2^{j+1}r}^{2^{j+2}r}\left(  1+\ln \frac{2^{j+1}r}{r}\right)
^{2}\left \vert B\left(  x_{0},2^{j+1}r\right)  \right \vert ^{\lambda
_{2}-\left(  \frac{1}{p_{1}}+\frac{1}{p_{2}}\right)  +2}\prod \limits_{i=1}%
^{2}\left \Vert f_{i}\right \Vert _{L_{p_{i}}\left(  B\left(  x_{0}%
,2^{j+1}r\right)  \right)  }\frac{dt}{t^{2n+1}}\\
& \lesssim \Vert b_{2}\Vert_{LC_{q_{2},\lambda_{2}}^{\left \{  x_{0}\right \}  }%
}\int \limits_{2r}^{\infty}\left(  1+\ln \frac{t}{r}\right)  ^{2}\left \vert
B\left(  x_{0},t\right)  \right \vert ^{\lambda_{2}-\left(  \frac{1}{p_{1}%
}+\frac{1}{p_{2}}\right)  +2}\prod \limits_{i=1}^{2}\left \Vert f_{i}\right \Vert
_{L_{p_{i}}\left(  B\left(  x_{0},t\right)  \right)  }\frac{dt}{t^{2n+1}}\\
& \lesssim \Vert b_{2}\Vert_{LC_{q_{2},\lambda_{2}}^{\left \{  x_{0}\right \}  }%
}\int \limits_{2r}^{\infty}\left(  1+\ln \frac{t}{r}\right)  ^{2}%
{\displaystyle \prod \limits_{i=1}^{2}}
\Vert f_{i}\Vert_{L_{p_{i}}(B(x_{0},t))}\frac{dt}{t^{n\left(  \left(  \frac
{1}{p_{1}}+\frac{1}{p_{2}}\right)  -\lambda_{2}\right)  +1}}.
\end{align*}
This completes the proof of inequality (\ref{f3}).

Now we turn to estimate $G_{42}$. Let $1<\tau<\infty$, such that $\frac{1}%
{p}=\frac{1}{q_{1}}+\frac{1}{\tau}$. Then, by H\"{o}lder's inequality and
(\ref{f3}), we obtain%
\begin{align*}
G_{42}  & =\left \Vert \left[  \left(  b_{1}-\left \{  b_{1}\right \}
_{B}\right)  \right]  T^{\left(  2\right)  }\left[  f_{1}^{\infty},\left(
b_{2}-\left \{  b_{2}\right \}  _{B}\right)  f_{2}^{\infty}\right]  \right \Vert
_{L_{p}\left(  B\left(  x_{0},r\right)  \right)  }\\
& \lesssim \left \Vert \left(  b_{1}-\left \{  b_{1}\right \}  _{B}\right)
\right \Vert _{L_{q_{1}}\left(  B\right)  }\left \Vert T^{\left(  2\right)  }
\left[  f_{1}^{\infty},\left(  b_{2}-\left \{  b_{2}\right \}  _{B}\right)
f_{2}^{\infty}\right]  \right \Vert _{L_{\tau}\left(  B\right)  }\\
& \lesssim%
{\displaystyle \prod \limits_{i=1}^{2}}
\Vert \overrightarrow{b}\Vert_{LC_{q_{i},\lambda_{i}}^{\left \{  x_{0}\right \}
}}r^{\frac{n}{p}}\int \limits_{2r}^{\infty}\left(  1+\ln \frac{t}{r}\right)
^{2}%
{\displaystyle \prod \limits_{i=1}^{2}}
\Vert f_{i}\Vert_{L_{p_{i}}(B(x_{0},t))}\frac{dt}{t^{n\left(  \left(  \frac
{1}{p_{1}}+\frac{1}{p_{2}}\right)  -\left(  \lambda_{1}+\lambda_{2}\right)
\right)  +1}}.
\end{align*}

Similarly, $G_{43}$ has the same estimate above, here we omit the details,
thus the inequality%
\begin{align*}
G_{43}  & =\left \Vert \left[  \left(  b_{2}-\left \{  b_{2}\right \}
_{B}\right)  \right]  T^{\left(  2\right)  }\left[  \left(  b_{1}-\left \{
b_{1}\right \}  _{B}\right)  f_{1}^{\infty},f_{2}^{\infty}\right]  \right \Vert
_{L_{p}\left(  B\left(  x_{0},r\right)  \right)  }\\
& \lesssim%
{\displaystyle \prod \limits_{i=1}^{2}}
\Vert \overrightarrow{b}\Vert_{LC_{q_{i},\lambda_{i}}^{\left \{  x_{0}\right \}
}}r^{\frac{n}{p}}\int \limits_{2r}^{\infty}\left(  1+\ln \frac{t}{r}\right)
^{2}%
{\displaystyle \prod \limits_{i=1}^{2}}
\Vert f_{i}\Vert_{L_{p_{i}}(B(x_{0},t))}\frac{dt}{t^{n\left(  \left(  \frac
{1}{p_{1}}+\frac{1}{p_{2}}\right)  -\left(  \lambda_{1}+\lambda_{2}\right)
\right)  +1}}.
\end{align*}
is valid.

Finally, to estimate $G_{44}$, similar to the estimate of (\ref{f3}), we have%
\begin{align*}
& \left \vert T^{\left(  2\right)  }\left[  \left(  b_{1}-\left \{
b_{1}\right \}  _{B}\right)  f_{1}^{\infty},\left(  b_{2}-\left \{
b_{2}\right \}  _{B}\right)  f_{2}^{\infty}\right]  \left(  x\right)
\right \vert \\
& \lesssim%
{\displaystyle \sum \limits_{j=1}^{\infty}}
\left(  2^{j}r\right)  ^{-2n}\left[
{\displaystyle \int \limits_{B\left(  x_{0},2^{j+1}r\right)  }}
\left \vert b_{1}\left(  y_{1}\right)  -\left \{  b_{1}\right \}  _{B}\right \vert
\left \vert f_{1}\left(  y_{1}\right)  \right \vert dy_{1}\right]  \left[
{\displaystyle \int \limits_{B\left(  x_{0},2^{j+1}r\right)  }}
\left \vert b_{2}\left(  y_{2}\right)  -\left \{  b_{2}\right \}  _{B}\right \vert
\left \vert f_{2}\left(  y_{2}\right)  \right \vert dy_{2}\right] \\
& \lesssim%
{\displaystyle \prod \limits_{i=1}^{2}}
\Vert \overrightarrow{b}\Vert_{LC_{q_{i},\lambda_{i}}^{\left \{  x_{0}\right \}
}}\int \limits_{2r}^{\infty}\left(  1+\ln \frac{t}{r}\right)  ^{2}%
{\displaystyle \prod \limits_{i=1}^{2}}
\Vert f_{i}\Vert_{L_{p_{i}}(B(x_{0},t))}\frac{dt}{t^{n\left(  \left(  \frac
{1}{p_{1}}+\frac{1}{p_{2}}\right)  -\left(  \lambda_{1}+\lambda_{2}\right)
\right)  +1}}.
\end{align*}
Thus, we have%
\begin{align*}
G_{44}  & =\left \Vert T^{\left(  2\right)  }\left[  \left(  b_{1}-\left \{
b_{1}\right \}  _{B}\right)  f_{1}^{\infty},\left(  b_{2}-\left \{
b_{2}\right \}  _{B}\right)  f_{2}^{\infty}\right]  \right \Vert _{L_{p}\left(
B\left(  x_{0},r\right)  \right)  }\\
& \lesssim%
{\displaystyle \prod \limits_{i=1}^{2}}
\Vert \overrightarrow{b}\Vert_{LC_{q_{i},\lambda_{i}}^{\left \{  x_{0}\right \}
}}r^{\frac{n}{p}}\int \limits_{2r}^{\infty}\left(  1+\ln \frac{t}{r}\right)
^{2}%
{\displaystyle \prod \limits_{i=1}^{2}}
\Vert f_{i}\Vert_{L_{p_{i}}(B(x_{0},t))}\frac{dt}{t^{n\left(  \left(  \frac
{1}{p_{1}}+\frac{1}{p_{2}}\right)  -\left(  \lambda_{1}+\lambda_{2}\right)
\right)  +1}}.
\end{align*}
By the estimates of $G_{4j}$ above, where $j=1$, $2$, $3$, $4$. We know that%
\begin{align*}
G_{4}  & =\left \Vert T_{\left(  b_{1},b_{2}\right)  }^{\left(  2\right)
}\left(  f_{1}^{\infty},f_{2}^{\infty}\right)  \right \Vert _{L_{p}\left(
B\left(  x_{0},r\right)  \right)  }\lesssim%
{\displaystyle \prod \limits_{i=1}^{2}}
\Vert \overrightarrow{b}\Vert_{LC_{q_{i},\lambda_{i}}^{\left \{  x_{0}\right \}
}}r^{\frac{n}{p}}\int \limits_{2r}^{\infty}\left(  1+\ln \frac{t}{r}\right)
^{2}%
{\displaystyle \prod \limits_{i=1}^{2}}
\Vert f_{i}\Vert_{L_{p_{i}}(B(x_{0},t))}\\
& \times \frac{dt}{t^{n\left(  \left(  \frac{1}{p_{1}}+\frac{1}{p_{2}}\right)
-\left(  \lambda_{1}+\lambda_{2}\right)  \right)  +1}}.
\end{align*}

Recalling (\ref{100}), and combining all the estimates for $G_{1},$ $G_{2}$,
$G_{3}$, $G_{4}$, we get%
\[
\left \Vert T_{\left(  b_{1},b_{2}\right)  }^{\left(  2\right)  }\left(
f_{1},f_{2}\right)  \right \Vert _{L_{p}\left(  B\left(  x_{0},r\right)
\right)  }\lesssim%
{\displaystyle \prod \limits_{i=1}^{2}}
\Vert \overrightarrow{b}\Vert_{LC_{q_{i},\lambda_{i}}^{\left \{  x_{0}\right \}
}}r^{\frac{n}{p}}\int \limits_{2r}^{\infty}\left(  1+\ln \frac{t}{r}\right)
^{2}%
{\displaystyle \prod \limits_{i=1}^{2}}
\Vert f_{i}\Vert_{L_{p_{i}}(B(x_{0},t))}\frac{dt}{t^{n\left(  \left(  \frac
{1}{p_{1}}+\frac{1}{p_{2}}\right)  -\left(  \lambda_{1}+\lambda_{2}\right)
\right)  +1}}.
\]
Therefore, Theorem \ref{Teo 5} is completely proved.
\end{proof}

\subsection{\textbf{Proof of Theorem \ref{teo15}.}}

\begin{proof}
To prove Theorem \ref{teo15}, we will use the following relationship between
essential supremum and essential infimum%
\begin{equation}
\left(  \operatorname*{essinf}\limits_{x\in E}f\left(  x\right)  \right)
^{-1}=\operatorname*{esssup}\limits_{x\in E}\frac{1}{f\left(  x\right)
},\label{5}%
\end{equation}
where $f$ is any real-valued nonnegative function and measurable on $E$ (see
\cite{Wheeden-Zygmund}, page 143). Indeed, since $\overrightarrow{f}\in
LM_{p_{1},\varphi_{1}}^{\{x_{0}\}}\times \cdots$ $\times$ $LM_{p_{m}%
,\varphi_{m}}^{\{x_{0}\}}$, by (\ref{5}) and the non-decreasing, with respect
to $t$, of the norm $%
{\displaystyle \prod \limits_{i=1}^{m}}
\Vert f_{i}\Vert_{L_{p_{i}}(B\left(  x_{0},t\right)  )}$, we get%
\begin{align}
&  \frac{%
{\displaystyle \prod \limits_{i=1}^{m}}
\Vert f_{i}\Vert_{L_{p_{i}}(B\left(  x_{0},t\right)  )}}%
{\operatorname*{essinf}\limits_{0<t<\tau<\infty}%
{\displaystyle \prod \limits_{i=1}^{m}}
\varphi_{i}(x_{0},\tau)\tau^{\frac{n}{p_{i}}}}\leq \operatorname*{esssup}%
\limits_{0<t<\tau<\infty}\frac{%
{\displaystyle \prod \limits_{i=1}^{m}}
\Vert f_{i}\Vert_{L_{p_{i}}(B\left(  x_{0},t\right)  )}}{%
{\displaystyle \prod \limits_{i=1}^{m}}
\varphi_{i}(x_{0},\tau)\tau^{\frac{n}{p_{i}}}}\nonumber \\
&  \leq \operatorname*{esssup}\limits_{0<\tau<\infty}\frac{%
{\displaystyle \prod \limits_{i=1}^{m}}
\Vert f_{i}\Vert_{L_{p_{i}}(B\left(  x_{0},t\right)  )}}{%
{\displaystyle \prod \limits_{i=1}^{m}}
\varphi_{i}(x_{0},\tau)\tau^{\frac{n}{p_{i}}}}\leq%
{\displaystyle \prod \limits_{i=1}^{m}}
\left \Vert f_{i}\right \Vert _{LM_{p_{i},\varphi_{i}}^{\{x_{0}\}}}.\label{9}%
\end{align}
For $1<p_{1},\ldots,p_{m}<\infty$, since $(\varphi_{1},\ldots,\varphi
_{m},\varphi)$ satisfies (\ref{50}) and by (\ref{9}), we have%
\begin{align}
&  \int \limits_{r}^{\infty}\left(  1+\ln \frac{t}{r}\right)  ^{m}%
{\displaystyle \prod \limits_{i=1}^{m}}
\Vert f_{i}\Vert_{L_{p_{i}}(B(x_{0},t))}\frac{dt}{t^{^{n\left(
{\displaystyle \sum \limits_{i=1}^{n}}
\frac{1}{p_{i}}-%
{\displaystyle \sum \limits_{i=1}^{n}}
\lambda_{i}\right)  +1}}}\nonumber \\
&  \leq \int \limits_{r}^{\infty}\left(  1+\ln \frac{t}{r}\right)  ^{m}\frac{%
{\displaystyle \prod \limits_{i=1}^{m}}
\Vert f_{i}\Vert_{L_{p_{i}}(B(x_{0},t))}}{\operatorname*{essinf}%
\limits_{t<\tau<\infty}%
{\displaystyle \prod \limits_{i=1}^{m}}
\varphi_{i}(x_{0},\tau)\tau^{\frac{n}{p_{i}}}}\frac{\operatorname*{essinf}%
\limits_{t<\tau<\infty}%
{\displaystyle \prod \limits_{i=1}^{m}}
\varphi_{i}(x_{0},\tau)\tau^{\frac{n}{p_{i}}}}{t^{n\left(
{\displaystyle \sum \limits_{i=1}^{n}}
\frac{1}{p_{i}}-%
{\displaystyle \sum \limits_{i=1}^{n}}
\lambda_{i}\right)  }}\frac{dt}{t}\nonumber \\
&  \leq C%
{\displaystyle \prod \limits_{i=1}^{m}}
\left \Vert f_{i}\right \Vert _{LM_{p_{i},\varphi_{i}}^{\{x_{0}\}}}%
\int \limits_{r}^{\infty}\left(  1+\ln \frac{t}{r}\right)  ^{m}\frac
{\operatorname*{essinf}\limits_{t<\tau<\infty}%
{\displaystyle \prod \limits_{i=1}^{m}}
\varphi_{i}(x_{0},\tau)\tau^{\frac{n}{p}}}{t^{n\left(
{\displaystyle \sum \limits_{i=1}^{n}}
\frac{1}{p_{i}}-%
{\displaystyle \sum \limits_{i=1}^{n}}
\lambda_{i}\right)  +1}}dt\nonumber \\
&  \leq C%
{\displaystyle \prod \limits_{i=1}^{m}}
\left \Vert f_{i}\right \Vert _{LM_{p_{i},\varphi_{i}}^{\{x_{0}\}}}\varphi
(x_{0},r).\label{13}%
\end{align}
Then by (\ref{200}) and (\ref{13}), we get%
\begin{align*}
\left \Vert T_{\overrightarrow{b}}^{\left(  m\right)  }\left(  \overrightarrow
{f}\right)  \right \Vert _{LM_{p,\varphi}^{\{x_{0}\}}} &  =\sup_{r>0}%
\varphi \left(  x_{0},r\right)  ^{-1}|B\left(  x_{0},r\right)  |^{-\frac{1}{p}%
}\left \Vert T_{\overrightarrow{b}}^{\left(  m\right)  }\left(  \overrightarrow
{f}\right)  \right \Vert _{L_{p}\left(  B\left(  x_{0},r\right)  \right)  }\\
&  \lesssim%
{\displaystyle \prod \limits_{i=1}^{m}}
\left \Vert \overrightarrow{b}\right \Vert _{LC_{q_{i},\lambda_{i}}^{\left \{
x_{0}\right \}  }}\sup_{r>0}\varphi \left(  x_{0},r\right)  ^{-1}\int
\limits_{r}^{\infty}\left(  1+\ln \frac{t}{r}\right)  ^{m}%
{\displaystyle \prod \limits_{i=1}^{m}}
\Vert f_{i}\Vert_{L_{p_{i}}(B(x_{0},t))}\frac{dt}{t^{^{n\left(
{\displaystyle \sum \limits_{i=1}^{n}}
\frac{1}{p_{i}}-%
{\displaystyle \sum \limits_{i=1}^{n}}
\lambda_{i}\right)  +1}}}\\
&  \lesssim%
{\displaystyle \prod \limits_{i=1}^{m}}
\left \Vert \overrightarrow{b}\right \Vert _{LC_{q_{i},\lambda_{i}}^{\left \{
x_{0}\right \}  }}%
{\displaystyle \prod \limits_{i=1}^{m}}
\left \Vert f_{i}\right \Vert _{LM_{p_{i},\varphi_{i}}^{\{x_{0}\}}}.
\end{align*}
Thus we obtain (\ref{51}). Hence the proof is completed.
\end{proof}

\end{document}